\documentclass{amsart}

\usepackage{graphicx}
\usepackage{subcaption}
\usepackage{hyperref}
\usepackage{xcolor}
\usepackage{seqsplit}
\usepackage{amsfonts}

\graphicspath{{figures/}}

\newtheorem{thm}{Theorem}
\newtheorem{lemma}{Lemma}

\title[Computation of Tight Enclosures for Laplacian Eigenvalues]
{Computation of Tight Enclosures for\\ Laplacian
  Eigenvalues}
\author{Joel Dahne\and Bruno Salvy}
\begin{document}
\maketitle
\begin{abstract}
  Recently, there has been interest in high-precision approximations
  of the first eigenvalue of the Laplace-Beltrami operator on
  spherical triangles for combinatorial purposes. We compute improved
  and certified enclosures to these eigenvalues. This is achieved by
  applying the method of particular solutions in high precision, the
  enclosure being obtained by a combination of interval arithmetic and
  Taylor models. The index of the eigenvalue is certified by
  exploiting the monotonicity of the eigenvalue with respect to the
  domain. The classically troublesome case of singular corners is
  handled by combining expansions at all corners and an expansion from
  an interior point. In particular, this allows us to compute~100
  digits of the fundamental eigenvalue for the 3D Kreweras model that
  has been the object of previous efforts.
\end{abstract}

\section*{Introduction}

The most classical situation for the computation of Laplacian
eigenvalues is that of the Laplacian on a bounded open set
$\Omega\subset\mathbb{R}^d$. There is an increasing sequence
$(\lambda_n)_n$ of positive real numbers, called the eigenvalues, and
corresponding (eigen)functions $(u_n)_n$ in $C^\infty(\Omega)$ such
that $\Delta u_n+\lambda_nu_n=0$ in $\Omega$, while $u_n=0$ on the
boundary $\partial\Omega$ (a Dirichlet condition). Moreover, $(u_n)$
is a Hilbert basis of $L^2(\Omega)$ (see, e.g.,
\cite{CourantHilbert1962,Brezis2011}). The method of particular
solutions was introduced by Fox, Henrici and
Moler~\cite{FoxHenriciMoler1967} and later refined by Betcke and
Trefethen~%
\cite{BetckeTrefethen2005}. Betcke~\cite{Betcke2008} describes this
method as ``especially effective for very accurate computations''. (He
also gives pointers to precursors and related works.) Starting from a
solution $u^*$ of $\Delta u+\lambda^*u=0$ in $\Omega$ for some
$\lambda^*$ that does not necessarily satisfy
$u^*|_{\partial\Omega}=0$, this method lets one deduce an interval
around $\lambda^*$ that contains an eigenvalue of the original problem
and whose diameter can be bounded in terms of
$\max_ {x\in\partial\Omega}|u^*(x)|$ and $\|u^*\|_2$. The candidate
pair $(\lambda^*,u^*)$ is computed by first finding a set
$(f_{n,\lambda})$ of solutions of $\Delta u+\lambda u=0$ in $\Omega$;
next, looking for a linear combination of unit norm which is minimal
on the boundary; finally looking for $\lambda$ where this minimum
value is close to~0. For the prototypical example of the L-shaped
region in the plane (displayed in Figure~\ref{fig:lshape-domain} on
p.~%
\pageref{fig:lshape-domain}), Fox, Henrici and Moler could compute 5
certified~digits of $\lambda_1$ with their method, using $10$ terms of
the linear combination and further tricks exploiting the symmetry of
the domain. With their improved method, Betcke and Trefethen produced
14~digits of $\lambda_1$ and certified 13 of them, without having to
exploit any special feature and using 60 terms of the linear
combination. In higher precision, this nice behaviour persists, as
already observed by Jones~\cite{Jones2017}. We show how certifying an
enclosure is also possible: these 14~digits can be certified with~100
terms of the expansion. The width of the enclosure scales well: using
180~terms gives 27 certified digits (see Section~\ref{sec:ex-L}).

A corner of a polygonal domain in the plane is called \emph{regular}
if it has an angle on
the form $\pi/k$ for some nonnegative integer $k$; corners with
angles not of that type are called \emph{singular}. The regularity of
a corner results in eigenfunctions that can be continued analytically
in a neighborhood of the corner by a reflection argument~%
\cite[V\S16.6]{CourantHilbert1962}. The L-shaped region in the plane
is an
example of a region with only one singular corner (the reentrant one).
This is a favorable situation for the method: the expansion at the
singular corner has no difficulty converging at the other ones. By
the same reflection argument, a similar phenomenon takes place in the
case of spherical triangles. This
lets us improve upon previous work in this setting for triangles with
at most one singular corner.

For polygonal domains with several singular corners, Betcke and
Trefethen~\cite{BetckeTrefethen2005} used expansions at all the
singular corners. In our experiments with spherical triangles, this
technique alone has not been sufficient to make the method converge.
However, there are known cases when, in different methods,
joining data from the corners with data from the interior of the
domain is successful~%
\cite{PlatteDriscoll2004,GopalTrefethen2019a}.
This also happens in our situation, where we have observed that taking
expansions at all singular
corners and complementing them with an expansion at an interior
point worked very well. That is how we could compute 100~digits for the
triangle with angles~$(2\pi/3,2\pi/3,2\pi/3)$, see
Section~\ref{sec:3dkreweras}.

\subsection*{Combinatorial motivation}
Recently, these computations have become relevant in the very active
study of discrete walks in~$\mathbb{N}^d$ (see recent
surveys~\cite{Krattenthaler2015a,Bostan2019} for numerous references).
The relation between the Brownian motion and the heat equation can be
exploited to derive the asymptotic number of walks in~$\mathbb{N}^d$
starting and ending at the origin and using $n$ steps, all taken from
a given finite set~$S\subset \mathbb{Z}^d$~\cite{DenisovWachtel2015}.
Under mild conditions, this number behaves asympotically like
\begin{equation}\label{eq:asymptfn}
f_S(n)\sim K\rho^nn^{\alpha},\quad \alpha=-1-\sqrt{\lambda_1+
(d/2-1)^2},
\end{equation}
where~$\lambda_1$ is the first eigenvalue (called the
\emph{fundamental eigenvalue}) of the Laplace-Beltrami operator on the
sphere~$\mathbb{S}^{d-1}$, with Dirichlet boundary conditions~0 on a
spherical cone that can be computed from the step set~$S$ (the
constant $\rho$, which is more important asymptotically, can also be
computed from~$S$). A question of interest in combinatorics is the
nature of the sequence~$f_S(n)$, i.e., the type of recurrence it may
satisfy. Depending on the step set~$S$, it can be solution of a linear
recurrence with constant coefficients, or with polynomial
coefficients, or of no such recurrence. The asymptotics above can be
used to rule out possibilities. For instance, if $f_S(n)$ satisfies a
linear recurrence with constant coefficients, then the
exponent~$\alpha$ has to be a nonnegative integer. A deeper result is
that if this integer sequence satisfies a linear recurrence with
\emph{polynomial} coefficients, then $\alpha$ has to be
\emph{rational}. This has been used to complete the classification of
\emph{planar} lattice walks with small steps (steps in
$\{0,\pm1\}^2$)~%
\cite{BostanRaschelSalvy2014}. It revealed a strong connection between
the existence of a linear recurrence with polynomial coefficients and
the finiteness of a group associated to the walk.

In dimension~3, an open question is whether the connection between
linear recurrences and finiteness of the associated group still holds.
Recently, Bogosel \emph{et alii} established that the cases where the
group is finite correspond to 17~triangles tiling the sphere that they
give explicitly~\cite{BogoselPerrollazRaschelTrotignon2020}. Of these,
7~triangles correspond to the small number of
spherical triangles with three regular corners for which the eigenvalues of
the Laplace-Beltrami operator are known explicitly~%
\cite{BerardBesson1980,Berard1983}. For the remaining 10~cases,
only
numerical approximations are available. Obviously, from a numerical
estimate of the fundamental eigenvalue~$\lambda_1$, one cannot expect
to obtain a guarantee of the rationality of the related
 exponent~$\alpha$ in Eq.~%
\eqref{eq:asymptfn}. Instead, our aim is to use
this numerical approximation as a filter by giving a lower bound on
the size the denominator would have if the exponent was a rational
number; small bounds would suggest the need for further combinatorial
investigation. Thus, this is a situation where we are interested in
computing tens or, ideally, hundreds of digits of the fundamental
eigenvalue of the Laplace operator \footnote{This is a possible answer
  to the conclusion of the review of Jones' article~\cite{Jones2017}
  on MathSciNet: \emph{What does one do with the thousands (or even
    hundreds) of digits for these eigenvalues?}}.

\subsection*{Certified computation}
Eigenvalues of self-adjoint operators and their computation form a
classical topic of numerical analysis; good surveys are available~%
\cite{KuttlerSigillito1984,Boffi2010,GrebenkovNguyen2013,Strohmaier2017}.

The possible influence of rounding errors in the computation of bounds
for eigenvalues leads naturally to the use of interval arithmetic in
certified computations. The development of methods that are most
suitable in this context started in the 1990's~%
\cite{Plum1991,BehnkeGoerisch1994,NakaoYamamotoNagatou1999}.
Since then, the theoretical aspects have been extended to more and
more general equations; see the recent book by Nakao, Plum and
Watanabe~\cite{NakaoPlumWatanabe2019} for an account of the main
methods. The case of the Laplace operator has been studied by Liu and
Oishi~\cite{LiuOishi2013}. They mention the method of particular
solutions, but discard it because it does not guarantee the index of
the eigenvalue. Instead, they construct explicit eigenvalue
bounds from a finite element method.

Our contribution to this problem
is to show how the method of particular solutions lends itself to
certified computations: we detail the extra work required to certify
an enclosure for the eigenvalues and their index and show that this is
not the time-consuming part of the computation. Moreover, that method
has the advantage that it can easily be used to obtain high-precision
results, which is difficult by methods based on finite elements
(see for instance the discussion at the end of~\cite{GopalTrefethen2019a}).

In order to compute an enclosure for the fundamental eigenvalue by the
method of particular solutions, the first step is to use the method
without worrying about certified computations: any approximate pair
$(\lambda^*,u^*)$ will do. The difficulties at this stage are the same
as in a classical computation, mainly the slow convergence in the
presence of singular corners and the linear combinations of basis
functions that are very close to~0 inside the domain and make the
linear algebra problem ill-conditioned. Once these difficulties are
overcome and such an approximation at high precision is obtained, and
only then, we need a more careful computation when bounding the
distance to an actual close-by eigenvalue so that the enclosure can be
guaranteed. This is described in Section~\ref{sec:rigorous}.

Finally, we also need to prove that the eigenvalue that has been
produced is indeed the fundamental one. For the spherical
triangles we study, this is done in
Section~\ref{sec:cert_index}, by exploiting the monotonicity property
of eigenvalues with respect to the domain, certified roots of Ferrers
functions and a variant of Sturm's theorem.

\subsection*{Recent works}
Three recent works are most directly related to ours.

Jones~\cite{Jones2017} uses the method of particular solutions in high
precision for the computation of eigenvalues for polygons in the
plane. In particular, he obtains 1\,000 digits of the fundamental
eigenvalue of the L-shaped region. The correct digits are obtained
thanks to an empirical observation: the eigenvalues of the truncated
problems obtained by taking $N$ points on the boundary alternate below
and above the limiting value as $N$ increases. In our context of
certified computation, we cannot rely on this heuristic approach. We
revert to computing certified bounds on the maximum value on the
boundary and the norm of approximate eigenfunctions.

Very recently, G\'omez-Serrano and Orriols \cite{gomez-serrano19:any}
have proved that three eigenvalues do not determine a triangle. For
this, they used the method of particular solutions in the plane in a
spirit very similar to ours. The main differences with our work is
that we need much higher precision, that the triangles they are
interested in do not have any
regular corner and that we work on the
sphere rather than on the plane. Also, they have to deal with many
more triangles. Their
way of lower bounding the norm of the candidate eigenfunction is
extended to singular spherical triangles in
Section~\ref{section:norm}.

The finite element method has the disadvantage that a large
discretization is needed in order to get a good accuracy. However,
together with extrapolation, this approach was used with success by
Bogosel, Perrollaz, Raschel and Trotignon~%
\cite{BogoselPerrollazRaschelTrotignon2020} in our problem. They
compute approximations of the fundamental eigenvalue for
all~17 spherical
triangles corresponding to walks with finite groups. Their results are
compared to ours in Tables~\ref{table:regular-triangles},
\ref{table:singular-triangles} below. We can vastly improve the
precision they obtain, both for triangles with at most one singular
corner and for more singular triangles.

\subsection*{Plan} In Section~\ref{sec:mps}, we first recall the
method of particular solutions. Then, in Section~\ref{sec:rigorous} we
spell out the steps we use in the certification stage and illustrate
these in detail in the classical case of the L-shaped region
(\S\ref{sec:ex-L}). This example is presented in such a way that the
same steps apply to the spherical triangles in
Section~\ref{sec:LBops}. We then show how the index of the eigenvalue
can be certified in Section~\ref{sec:cert_index}. The results for the
spherical triangles that had been considered by Bogosel et al. are
then given (\S\ref{sec:regular},\S\ref{sec:singular-triangles}) and in
\S\ref{sec:denominators}, we conclude by providing lower bounds on the
denominators the corresponding exponents from Eq.~\eqref{eq:asymptfn}
would have if they were rational numbers, the case of the 3D Kreweras
model being detailed in Section~\ref{sec:3dkreweras}.

\section{Method of Particular Solutions}
\label{sec:mps}
The starting point of the method is the following a posteriori bound.
It allows one to enclose an eigenvalue by finding good approximations
to the eigenfunction.
\begin{thm}
  \label{thm:FoxHenriciMoler}
  \cite{FoxHenriciMoler1967,MolerPayne1968} Let
  $\Omega\subset\mathbb{R}^n$ be bounded. Let $\lambda$ and $u$ be
  an approximate eigenvalue and eigenfunction---that is,  they
  satisfy
  $\Delta u+\lambda u=0$ in $\Omega$ but not necessarily $u=0$
  on~$\partial\Omega$. Define
  \begin{equation*}
    \epsilon = \frac{\sqrt{\operatorname{Vol}(\Omega)}\sup_{x \in
    \partial\Omega}|u(x)|}{\|u\|_2}.
  \end{equation*}
  Then there exists an eigenvalue $\lambda_*$ such that
  \begin{equation}
    \label{eq:eigenvalue-bound}
    \frac{|\lambda - \lambda_*|}{\lambda_*} \leq \epsilon.
  \end{equation}
\end{thm}
Variants of this result hold more generally for elliptic differential
operators. As stated here, the bound is due to Moler and Payne~%
\cite{MolerPayne1968}, improving upon the original bound by Fox,
Henrici and Moler~\cite{FoxHenriciMoler1967}. Further generalizations
and improvements to this bound are available in the
literature~\cite{KuttlerSigillito1984,Still1988,BarnettHassell2011,Strohmaier2017}.

\smallskip

The version of the method of particular solutions we give here is due
Betcke and Trefethen \cite{BetckeTrefethen2005,Betcke2008}. The
starting point is a set of solutions $u_{\lambda}^{(k)}$ to
$\Delta u + \lambda u = 0$ in~$\Omega$, that are not constrained
on the boundary~$\partial\Omega$. An approximate
eigenfunction is given by a linear combination of these,
\begin{equation*}
  u_{\lambda}^{*} = \sum_{k = 1}^{N}c_{k}u_{\lambda}^{(k)},
\end{equation*}
where the coefficients are chosen so that it is close to unit norm in
$\Omega$ and minimal on the boundary.

The coefficients $c_{k}$ are found by taking $m_{B}$ points on the
boundary, $\{x_{i}\}_{i = 1}^{m_{B}} \subset \partial \Omega$, on
which we want to minimize $u_\lambda^*$. This alone is not sufficient:
increasing the number $N$ of elements in the basis leads to the
existence of linear combinations very close to~0 inside the
domain~$\Omega$. These linear combinations lead to spurious solutions,
even away from an eigenvalue. Betcke and Trefethen cure this problem
by also considering $m_{I}$ points in the interior
$\{y_{j}\}_{j = 1}^{m_{I}} \subset \Omega$ and using them to make
sure that the linear combinations of interest stay close to unit norm
in the domain. The computation now involves two matrices
\begin{align*}
  (A_B)_{ik}(\lambda) &= u^{(k)}_\lambda(x_i), &&1\le i\le m_B,\ 1\le
  k\le N;\\
  (A_I)_{jk}(\lambda) &= u^{(k)}_\lambda(y_j), &&1\le j\le
                                                  m_I,\ 1\le k\le N.
\end{align*}

The matrices $A_B$ and $A_I$ can be combined into a matrix whose
$QR$ factorization gives an orthonormal basis of these function
evaluations
\begin{equation}\label{eq:QR}
  A(\lambda) = \begin{bmatrix} A_B(\lambda)\\ A_I(\lambda) \end{bmatrix}
= \begin{bmatrix} Q_B(\lambda)\\ Q_I(\lambda) \end{bmatrix}R
(\lambda)=:Q(\lambda)R(\lambda).
\end{equation}
Now, the right singular vector~$v$ of norm~1 corresponding to the
smallest singular value $\sigma(\lambda)$ of the top part
$Q_{B} (\lambda)$ is a good candidate for the eigenfunction when
$\sigma(\lambda)$ is small. Moreover,
\[\|Q(\lambda)v\|_2^2=1=\sigma(\lambda)^2+\|Q_I(\lambda)v\|^2,\]
forcing the function to have norm close to~1 on the interior points.
Thus the next step of the method is to search
for $\lambda$ minimizing~$\sigma(\lambda)$. Once such a $\lambda$
is found, the actual vector $c$ of coefficients is
recovered by solving the linear system
$ R(\lambda)c = v$. As noted by Betcke and Trefethen this linear
system is highly ill-conditioned but the solutions obtained are
nevertheless small on the boundary of the domain. We observe the same
phenomenon in our computations, taking place at a lot higher
precision. Computing the minimum to a certain precision only requires
10-20 extra bits of precision in the intermediate computations even
for as high as 300 bits of precision for the minimum.

Furthermore, Betcke and Trefethen give a geometric interpretation
of~$\sigma(\lambda)$: it is the sine of the angle between the subspace
of $\mathbb{R}^{m_B+m_I}$ generated by the columns of~$A$~---~the
values of the functions~$ (u_\lambda^{(k)})$ at the interior and
boundary points~$x_i$ and~$y_i$~---~and the
subspace~$ \mathbb{R}^ {m_B}\times 0^{m_I}$~---~the values of
functions that are~0 at the boundary points~$y_i$. This angle
becomes~0 when these spaces intersect, which means that a combination
of the functions~$ (u_\lambda^{(k)})$ is~0 at the boundary
points~$y_i$. Betcke shows that the tangent of that same angle is the
smallest generalized singular eigenvalue of the pencil~$\{A_B,A_I\}$,
which gives another way of computing it~\cite{Betcke2008}.

\section{Certified computation of the enclosure}
\label{sec:rigorous}
The method of particular solutions can be performed in high precision
environments as provided by computer algebra systems or by specialized
libraries like MPFR, but of course, running in high precision does not
guarantee anything about the accuracy of the result. Given an
approximate
eigenfunction $u$, a certified enclosure of the eigenvalue is
provided by Theorem~\ref{thm:FoxHenriciMoler} provided we obtain an
upper bound on the maximum of $|u|$ on the boundary
$\partial\Omega$ and a lower bound on its norm in $\Omega$.

As in the previous section, we let $u(x)$ be given by a linear
combination of functions $u_{k}(x)$
\begin{equation*}
  u(x) = \sum_{k = 1}^{N}c_{k}u_{k}(x).
\end{equation*}
The details of the implementation of the method depend on the actual
family of functions $u_{k}(x)$ and on the domain. The following
approach works in all our examples that are two-dimensional.

\subsection{Upper bound on the boundary}
\label{section:upper-bound}
This is the difficult part. A key role in certified computations is
played by interval arithmetic~\cite{Tucker2011,Johansson2017}, but it
cannot be applied blindly. The function $u(x)$ is given by a sum of
individually large terms whose sum is expected to be very small on the
boundary. A direct use of interval arithmetic is bound to fail, as it
handles cancellations poorly. The technique of Taylor
models~\cite{MakinoBerz2003} overcomes this difficulty: instead of
bounding the function directly, one computes a Taylor expansion of the
function before computing a bound by interval arithmetic. If
$\gamma: [0,1]\rightarrow\partial\Omega$ is a parametrization of the
boundary and $I\subset[0,1]$ is an interval of width~$|I|$, then
\begin{equation*}
  \max_{t\in I}u(\gamma(t)) \leq \max_{t \in I} P_{\ell - 1}(t)
  + \frac{(|I|/2)^{\ell}}{\ell!}\max_{t \in I}
  \left|\frac{d^{\ell}}{dt^{\ell}}u(\gamma(t))\right|,
\end{equation*}
where $P_{\ell - 1}$ is the Taylor polynomial of $u(\gamma(t))$ of
order $\ell - 1$ centered around the midpoint of the interval. Since
the functions~$u_k$ are obtained from the Laplacian by separation of
variables, they satisfy linear differential equations from which
linear recurrences for the Taylor coefficients follow. Thus the
polynomial~$P_{\ell - 1}$ can be computed efficiently to high
precision. Bounding it on the interval $I$ can be done by locating
possible extrema using classical interval arithmetic. Bounding
$\left|\frac{d^{\ell}}{dt^{\ell}}u(\gamma(t))\right|$ on the
interval~$I$ by interval arithmetic does suffer from the same problems
as bounding $u(\gamma(t))$ directly. This is however mitigated by the
factor ${(|I|/2)}^{\ell}/{\ell!}$, which can be made small by choosing
a sufficiently large~$\ell$ and/or a sufficiently small~$|I|$.

Without any extra computational effort, this also gives a lower bound
on the boundary as
\begin{equation*}
  \min_{t \in I} u(\gamma(t)) \geq \min_{t \in I} P_{\ell - 1}(t)
  - \frac{(|I|/2)^{\ell}}{\ell!}\max_{t \in I}
  \left|\frac{d^{\ell}}{dt^{\ell}}u(\gamma(t))\right|,
\end{equation*}
which can be used when lower bounding the norm below.

In the computations, $\ell$ is chosen to be the same as the number of
terms in the expansion of $u$. This choice is heuristic but in
practice this means that the number of subdivisions required stays
approximately constant as the number of terms in $u$ increases. Due to
the fact that there are a lot of cancellations, using sufficiently
high precision in the computations is important. In practice we add 30
bits of precision for computations with a target precision of 100-300
bits, only the high precision computations of the Kreweras triangle in
Section~\ref{sec:3dkreweras} require up to 80 extra bits of
precision. Adding more precision is however rather cheap and we have
favored using unnecessary bits over trying to minimize the number of
bits used in intermediate computations.

\subsection{Lower bound on the norm}\label{section:norm}
This is easier as only the magnitude is needed and no cancellation
takes place. The norm can be lower bounded by computing it on a subset
of the domain, $\Omega' \subset \Omega$. In our examples two
different cases occur.

In the case when there is a single expansion in an orthogonal basis,
it is natural to let $\Omega'$ be a circular sector inside the
domain. Then the integral giving the norm of $u(x)$ splits into a
sum of one-dimensional integrals. These one-dimensional integrals can
then be lower bounded by means of a verified
integrator~\cite{JohanssonMezzarobba2018}. This is what we do for the
L-shaped region below and in the case of regular triangles in
Section~\ref{sec:regtriglowerbound}.

When the terms in the expansion are not orthogonal, more work is
required. We apply a method used by G\'omez-Serrano and Orriols
\cite{gomez-serrano19:any} in the context of polygons. The idea is
that if $u$ does not vanish in a subset $\Omega'$ of $\Omega$, without
loss of generality it can be assumed to be positive there and then,
since $-\Delta u =\lambda u>0$, $u$ is a superharmonic function and
satisfies $\inf_{\Omega'}u\ge\inf_{\partial\Omega'}u$. Thus in that
situation, a lower bound for $|u|$ on $\partial\Omega'$ yields a lower
bound for $|u|$ inside $\Omega'$.

In order to detect that $u$ does not vanish inside $\Omega'$, the
first step is to compute $\min u$ and $\max u$ on $\partial\Omega'$ as
above. If that does not allow one to decide that the sign of $u$ is
fixed
on this boundary, then~0 is the only lower bound we can deduce for
$|u|$ in $\Omega'$. Since only a lower bound is desired, this
technique is applied to subdomains of $\Omega$ that avoid its
boundaries and, if more precision is necessary, to subdivisions of
them.

The remaining task is to ensure that $u$ has fixed sign in
$\Omega'$ when it has a fixed sign on its boundary, which we assume to
be positive without loss of generality. The key observation is that
$u$ cannot be negative inside $\Omega'$ if $\Omega'$ is small enough.
Indeed, if $\Omega''\subset\Omega'$ is a maximal domain where $u<0$,
then $u=0$ on $\partial\Omega''$ and thus $\lambda$ is an eigenvalue
for $\Omega''$. In $\mathbb{R}^n$, the Faber-Krahn inequality states
that the ball minimizes the first Dirichlet eigenvalue among all
domains of the same volume. The situation for domains on the sphere is
similar~ \cite{AshbaughLevine1997}:
\begin{equation}
  \label{eq:faber-krahn}
  \lambda=\lambda_1(\Omega'')\ge\lambda_1(\Omega^\star),
\end{equation}
where $\Omega^\star$ is the spherical cap with the same area as
$\Omega''$. This eigenvalue can be computed explicitly using zeros of
the Legendre functions. Since it increases when the domain decreases,
for too small a domain~$\Omega'$ it cannot be smaller than $\lambda$.
Given~$\lambda$, one can precompute once and for all the size of the
regions~$\Omega'$ over which this method is not sufficient to
reach a conclusion, and subdivide those into smaller regions.

\subsection{Certifying the index of the eigenvalue} At this stage,
Theorem~\ref{thm:FoxHenriciMoler} asserts that an actual eigenvalue
lies at a controlled distance from the original approximation. It is
also possible to certify that this eigenvalue is indeed the
fundamental one. Since all the eigenvalues monotonically decrease when
the domain is enlarged it is sufficient to find a larger domain
containing~$\Omega$ for which the second eigenvalue is easy to compute
and is larger than our estimate. This part of the computation depends
on the domain, we refer to the sections below for instantiations of
this method in our examples.

\subsection{Notes about implementation}
\label{sec:implementation}
Our code\footnote{Available at
  \url{https://github.com/Joel-Dahne/MethodOfParticularSolutions.jl}.\\
  The implementation is single threaded and in all cases where timings
  are given the computations have been done on a relatively old Intel
  Xeon E5-2620 running at 2GHz.} is implemented in
Julia~\cite{Julia-2017} and relies on Arb~\cite{Johansson2017}, used
through its Julia interface in Nemo~\cite{nemo}, for most of the
numerics. Arb uses arbitrary precision ball arithmetic, which allows us
to work with very high precision and also certify the errors in the
computations. Its verified integrator allows us to compute a lower
bound of the norm. It also has support for computations with
polynomials and implements Taylor expansions for many common
functions. Some of the functions for which we require Taylor
expansions are not implemented; these are the Bessel functions and the
Ferrers functions. These functions do however satisfy linear
differential equations from which linear recurrences for their Taylor
coefficients can easily be obtained with the Maple package
{Gfun}~\cite{SalvyZimmermann1994}, which we used to implement them on
top of Arb, letting us compute Taylor polynomials of arbitrarily large
degree for the required functions. All this enables us to compute a
certified upper bound of $\epsilon$ from
Theorem~\ref{thm:FoxHenriciMoler}.

When computing the approximate eigenfunction using the method of
particular solutions we do not require certified computations. However
Arb is one of the few libraries which implements efficient arbitrary
precision versions of the special functions we require and we
therefore use it also to compute the matrix $A(\lambda)$ from
Equation \eqref{eq:QR}. The QR factorization and SVD are computed
directly in Julia, which relies on
MPFR~\cite{FousseHanrotLefevrePelissierZimmermann2007} for its
BigFloat type. The minimum of $\sigma(\lambda)$ is then found using
an implementation of Brent's method in Julia~\cite{mogensen2018optim}.

The program works by finding the minimum of $\sigma(\lambda)$ for a
fixed number of terms in the expansion and then iteratively increasing
the number of terms, giving better and better approximations of
$\lambda$. In principle we only need to compute the enclosure in the
final step, once we have found what we believe to be a sufficiently
good approximation. This means that the computational cost is the same
as for the regular method of particular solutions, plus the extra cost
of computing the enclosure in the end. However in the examples below
we choose to compute the enclosure at every step to make it easier to
follow the progress.

\subsection{Example of an L-shaped region}\label{sec:ex-L}
\begin{figure}
  \centering
  \includegraphics[height=4cm]{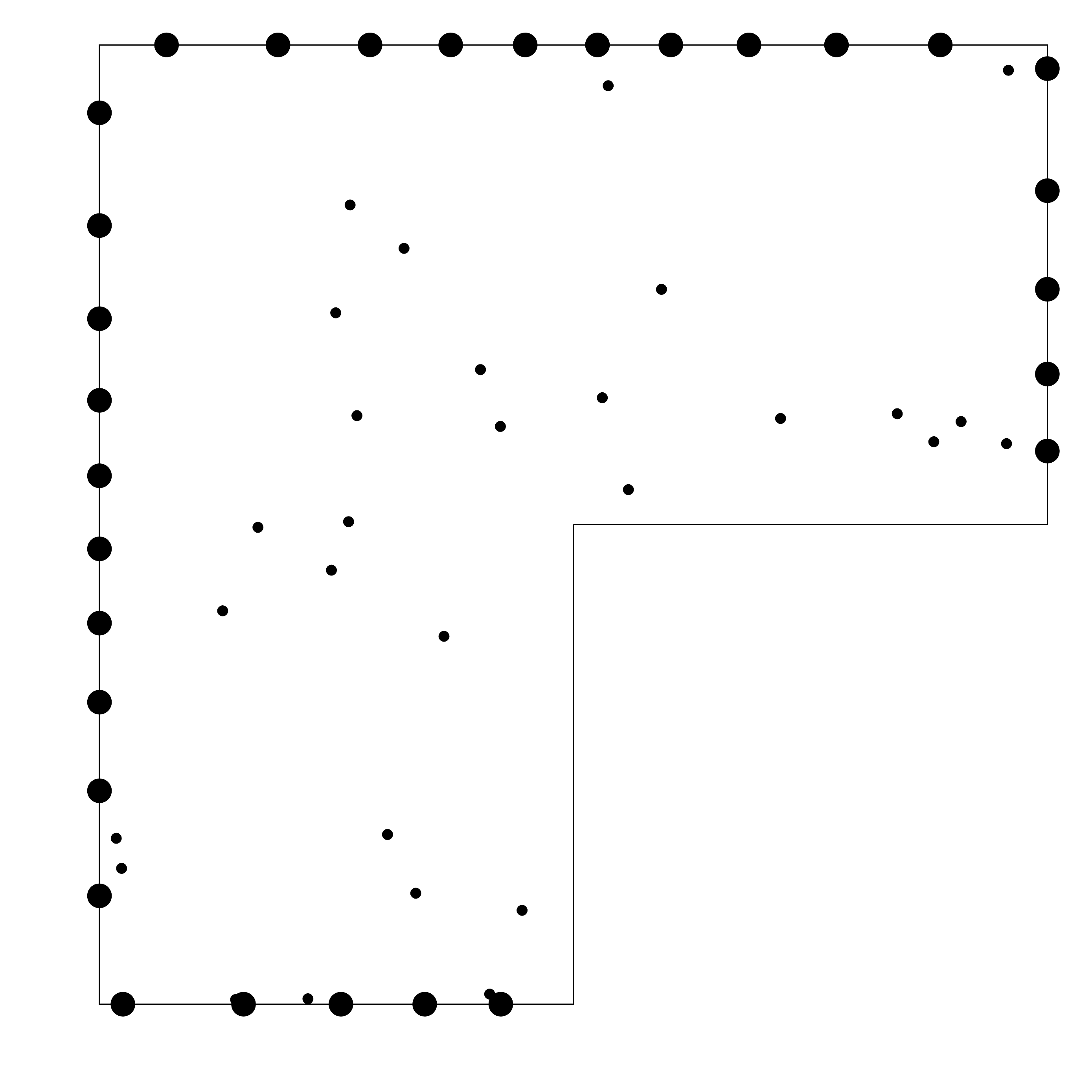}
  \caption{The L-shaped region with points on the boundary as well as
    random points the interior.}
  \label{fig:lshape-domain}
\end{figure}

The L-shaped region presented in Figure~\ref{fig:lshape-domain} is the
union of three squares of unit size. It is a classical test for
methods computing eigenvalues of the Laplacian~%
\cite{FoxHenriciMoler1967,DesclouxTolley1983,Still1988,PlatteDriscoll2004,BetckeTrefethen2005,Betcke2008,LiuOishi2013,Liu2015,Jones2017,CancesDussonMadayStammVohralik2017}.
We use it to exemplify the method of particular solutions on this
domain and the steps needed to give a certified enclosure, before
turning to spherical triangles in the next section. The presentation
is intentionally similar to that given by Betcke and Trefethen
\cite{BetckeTrefethen2005} so that the certification steps can be seen
clearly.

Solutions of $\Delta u+\lambda u=0$ in the plane are found by the
method of separation of variables: if $u=f(r)g(\theta)$ in polar
coordinates, then
\[r^2\frac{\Delta u}u=\frac{g''(\theta)}{g(\theta)}+r^2\frac{f''(r)}
  {f(r)}+r\frac{f'(r)}{f(r)},\] so that the equation
$\Delta u+\lambda u =0$ decouples into
\begin{equation}\label{eq:Bessel}
  g''(\theta)+K g(\theta)=0,\quad
  r^2f''(r)+rf'(r)+(\lambda
  r^2-K)f
  (r)=0
\end{equation}
for an arbitrary constant~$K$. The L-shaped domain has an angle
$3\pi/2$ at its reentrant corner, set at the origin. The boundary
conditions can be chosen so that the solutions are identically equal
to zero on the adjacent line segments and also finite at the origin.
Then, the first condition forces $K=4k^2/9$ with $k\in\mathbb Z$. The
second equation of Eq.~\eqref{eq:Bessel} is a variant of Bessel's
equation~%
\cite[Eq.~10.2.1]{OlverLozierBoisvertClark2010}, and therefore the
solutions that are also finite at the origin are given by
\[g(\theta)=\sin(2k\theta/3),\quad f(r)=J_{2k/3}(\sqrt{\lambda}r),\]
where $J_{2k/3}$ is the Bessel function of the first kind and
$k\in\mathbb{N}$.
With this basis, the linear combination
\begin{equation}\label{eq:expansion-u}
  u(r, \theta) = \sum_{k = 1}^{N}c_{k}u_{k}(r, \theta),\qquad
  \text{with}\quad
  u_{k}(r, \theta) = \sin(2k\theta/3)J_{2k/3}(\sqrt{\lambda}r)
\end{equation}
only has to be minimized on the remaining four boundary segments.

\subsubsection{Computation of a candidate}
\begin{figure}
  \centering
  \begin{subfigure}[t]{0.45\textwidth}
    \includegraphics[width=\textwidth]{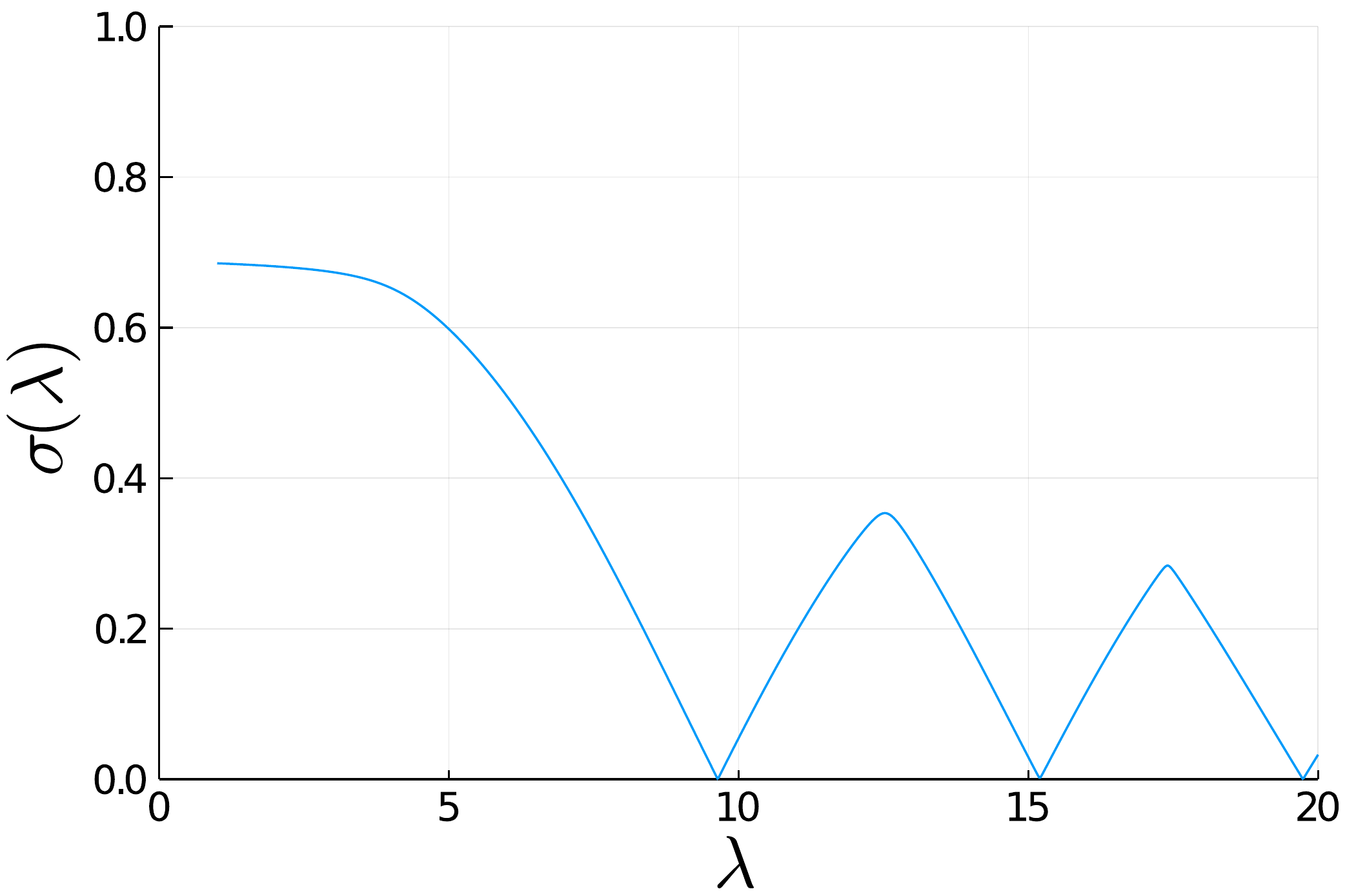}
    \caption{The function $\sigma(\lambda)$.}
    \label{fig:lshape-sigma}
  \end{subfigure}
  \hspace{0.05\textwidth}
  \begin{subfigure}[t]{0.45\textwidth}
    \includegraphics[width=\textwidth]{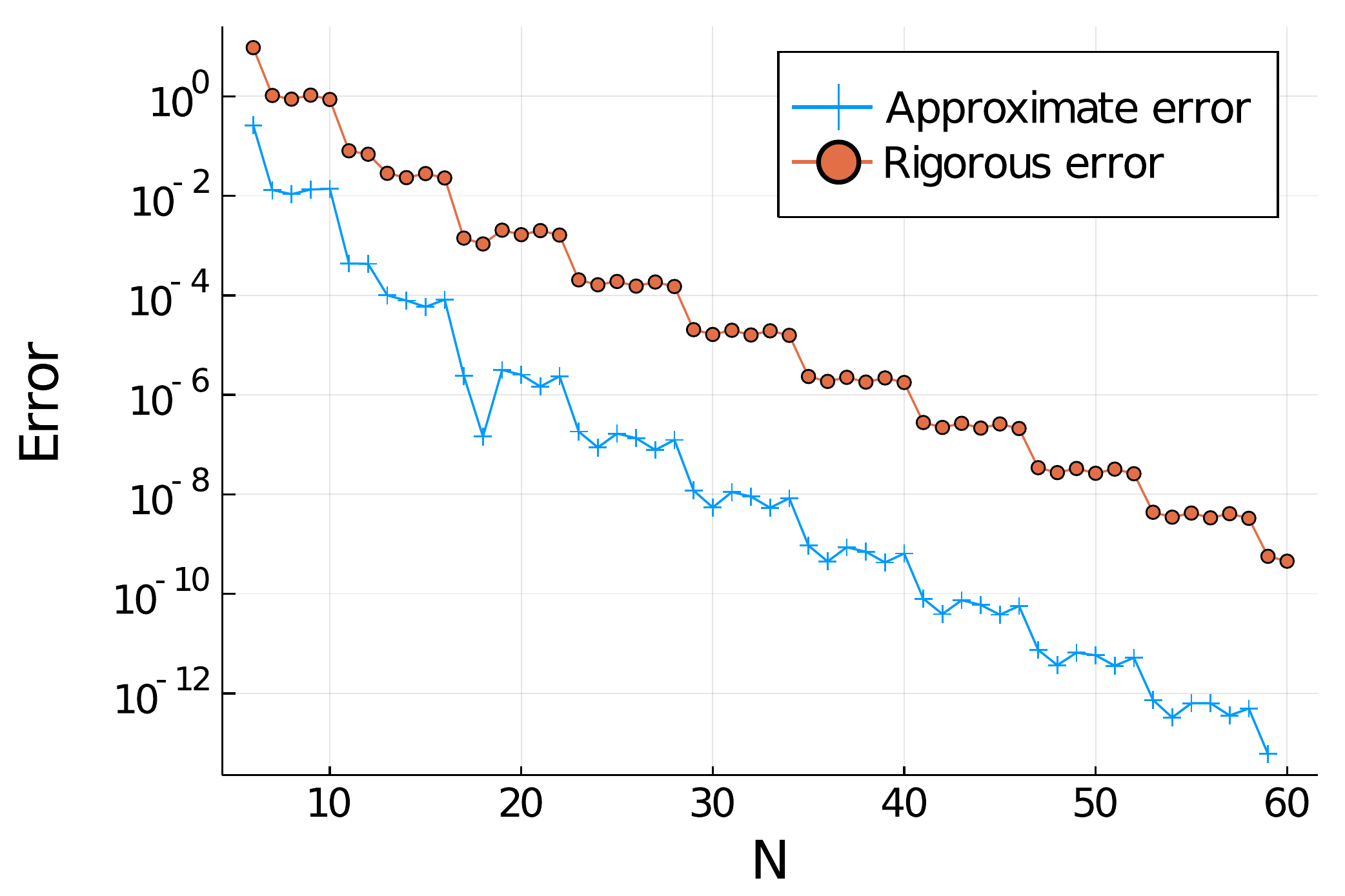}
    \caption{Convergence towards the fundamental eigenvalue.}
    \label{fig:lshape-convergence}
  \end{subfigure}
  \caption{Results for the L-shaped domain.}
\end{figure}

With the notation of Section~\ref{sec:mps}, we take $m_i=32$ random
points in the interior of the domain and $m_b=32$ points on the
boundary skipping the sides next to the reentrant corner, as shown in
Figure~\ref{fig:lshape-domain}. For values of $\lambda$ in the
interval~$[0,20]$, the functions $u_k$ from
Equation~\eqref{eq:expansion-u} are evaluated at those points and the
smallest singular value~$\sigma (\lambda)$ of the top part of the $QR$
factorization is computed. The resulting graph when using~$N=16$ terms
of the expansion is shown in Figure~\ref{fig:lshape-sigma}. In the
graph one can see the three minima corresponding to the first three
eigenvalues.

\subsubsection{Upper bound on the boundary}
The upper bound on the boundary is computed as described in
Section~\ref{section:upper-bound}, using truncated Taylor expansions
with certified bounds on the remainders.

\subsubsection{Lower bounding the norm}
We lower bound the norm of $u$ by considering the disk sector $G$
with radius 1 and angle $3\pi/2$ inscribed in the domain:
\begin{equation*}
  \|u\|^{2} =\int_\Omega{u^2\,dx}\ge \int_{G}u^{2} dx
  = \int_{0}^{3\pi/2}\int_{0}^{1}ru(r, \theta)^{2}dr d\theta.
\end{equation*}
When $u$ is given as a linear expansion of the form~%
\eqref{eq:expansion-u},
orthogonality of the family $\sin(2k\theta/3)$ on the interval
$ [0,3\pi/2]$ simplifies this last integral to
\begin{align*}
  \|u\|^{2}
  &\ge \sum_{k = 1}^{N}c_{k}^{2}\int_{0}^{3\pi/2}{\!\!\!\sin
    (2k\theta/3)^{2}d\theta}
    \int_{0}^{1}{rJ_{2k/3}(\sqrt{\lambda}r)^{2} dr}\\
  &= \frac{3\pi}{4}\sum_{k = 1}^{N}c_{k}^{2}\int_{0}^
    {1}{rJ_{2k/3}(\sqrt{\lambda}r)^{2} dr}.
\end{align*}
The remaining integrals
$\int_{0}^{1}rJ_{2k/3}(\sqrt{\lambda}r)^{2} dr$ can now be
efficiently computed with a certified integrator, giving a lower bound
for the norm. In practice we compute the integral from $\epsilon$ to
1 for a small $\epsilon$ to avoid having to deal with the branch cut
at 0.

\subsubsection{Convergence}
\begin{figure}
  \centering
  \includegraphics[height=5cm]{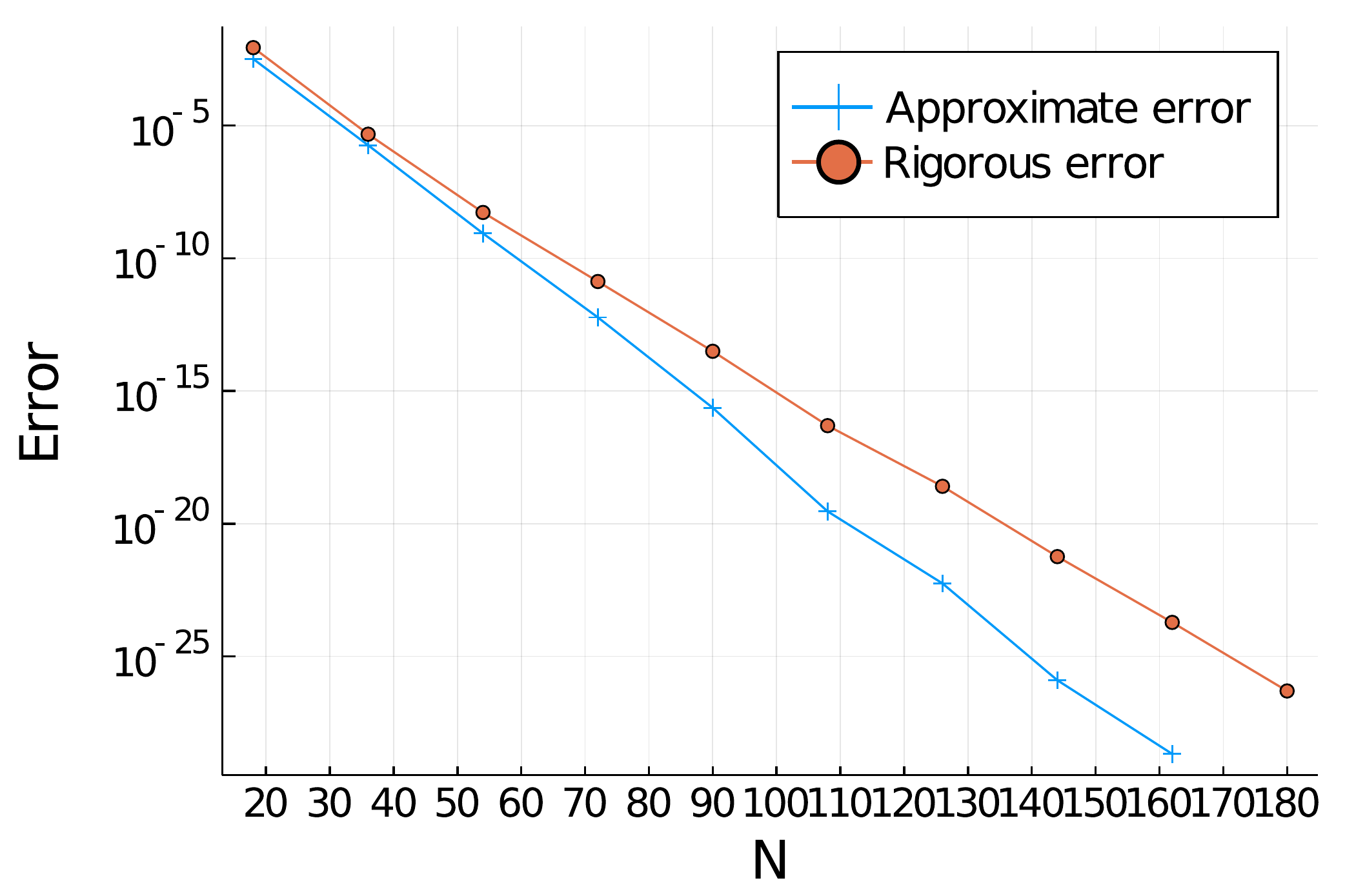}
  \caption{High precision results for convergence towards the
    fundamental eigenvalue for the L-shaped domain.}
  \label{fig:lshape-high-precision}
\end{figure}
The convergence towards the first minimum~---~the fundamental
eigenvalue~---~ is shown in Figure~\ref{fig:lshape-convergence}. It
shows two kinds of convergence. The first one is the approximate error
computed as the difference to the ``exact'' solution obtained with
$N = 60$. The other one is the radius of the certified enclosure,
giving an upper bound on the distance to the exact eigenvalue. While
the approximate error wobbles, the certified error is much more stable.
We also note that the quotient between the approximate and the
certified error remains mostly constant, the precision ``lost'' when
going from the approximate error to the certified error is about 3-4
digits for all $N$ greater than 25. This means that as we compute more
and more digits, the relative cost of considering the certified
enclosure instead of the approximate error decreases. For $N = 60$
we get the enclosure $\lambda \in [9.639723844 \pm 4.76\,10^{-10}]$.

As Betcke and Trefethen note, the problem of determining the
coefficients~$c_k$ is highly ill-conditioned. For $N = 60$ the
condition
number of $R(\lambda)$ in Eq.~\eqref{eq:QR} is about $10^{37}$. The values of
$A(\lambda)c$
can differ significantly from $Q(\lambda)\tilde{v}$ but are in general
still small. This ill-conditioning is not a problem for the
computation of the enclosure, since it is not important that we have
the
``correct'' solution, only that the computed solution is small on the
boundary.

Since the computations are done with arbitrary precision arithmetic we
can go further and compute more digits.
Figure~\ref{fig:lshape-high-precision} shows the convergence for $N$
up to 180. To avoid having to compute the enclosure many times we
increase $N$ in steps of 18; other than that the method is the same as
in Figure~\ref{fig:lshape-convergence}. These computations take
longer: 58 minutes, of which 48\% of the time is spent computing the
enclosure and the rest on finding the minimum of $\sigma(\lambda)$.
The final certified enclosure is
\[\lambda \in [9.63972384402194105271145926 \pm 7.36\,10^{-27}].\]

\section{Laplace-Beltrami Operators on Spherical Triangles}
\label{sec:LBops}
We now follow the same steps as for the L-shaped region for the
Laplace-Beltrami operator on the sphere.
\subsection{Laplace-Beltrami Operator on the Sphere}
In~$\mathbb R^d$, the Laplacian can be written
\[\Delta=r^{1-d}\frac\partial{\partial r}r^{d-1}\frac{\partial}
{\partial r}+r^{-2}\Delta_{\mathbb S^{d-1}},\]
where $\Delta_{\mathbb S^{d-1}}$ is the Laplace-Beltrami operator on
the $d$-dimensional sphere. In dimension~$d=3$ and spherical
coordinates, with $\theta$ denoting the polar angle and $\phi$ the
azimuthal angle, it becomes
\[\Delta_{\mathbb{S}^2} f(\theta,\phi)=\frac1{\sin\theta}
\frac{\partial}{\partial\theta}\!\left(\sin\theta\frac{\partial f}
{\partial\theta}\right)+\frac{1}{\sin^2\theta}\frac{\partial^2f}{\partial\phi^2}.\]
This is the operator whose fundamental eigenvalue on spherical
triangles is of interest for the asymptotics of lattice walks.

\subsection{Basis of Solutions}
It is classical that separation of variables applies: if
$u (\theta,\phi)=f(\theta)g (\phi)$, then
\[\sin^2(\theta)\frac{\Delta_{\mathbb S^{d-1}}u}u=\frac{g''(\phi)}{g
(\phi)}+\sin^2(\theta)\frac{f''(\theta)}{f
(\theta)}+\sin\theta\cos\theta\frac{f'(\theta)}{f(\theta)}.\]
Thus the equation $\Delta_{\mathbb S^{d-1}}u+\lambda u=0$ decouples
into
\begin{equation}\label{eq:legendre}
  g''(\phi)+K g(\phi)=0,\quad
  (1-x^2)f''(x)-2xf'(x)+\left(\lambda-\frac K{1-x^2}\right)f(x)=0,
\end{equation}
where $\cos\theta=x$, for an arbitrary constant~$K$. This is the
analogue of Equation~\eqref{eq:Bessel} in the previous section.

Let $D$ be a spherical triangle. Choose one vertex to place at the
north pole. Fix one of the sides on the meridian~$\phi=0$ and the
other one on the meridian~$\phi=\pi/\alpha$. Selecting solutions that
are identically~0 on both these meridians fixes~$K=\alpha^2k^2$ with
$k\in\mathbb Z$. The second equation in Eq.~\eqref{eq:legendre} is the
associated Legendre
equation~\cite[Eq.~14.2.2]{OlverLozierBoisvertClark2010}. Its
solutions that are real on the interval $(-1,1)$ and finite at $x=1$
(corresponding to the north pole) are the \emph{Ferrers functions of
  the first kind} $\mathsf{P}^{\mu}_{\nu}(x)$ with $\mu=-\alpha k$,
$k\in\mathbb N$ and $\nu$ obeying $\lambda=\nu(\nu+1)$~%
\cite[Eq.~14.8.1]{OlverLozierBoisvertClark2010}. Thus we are looking
for coefficients~$c_k$ such that the sum
\begin{equation}\label{eq:expansion-Ferrers}
  u(\theta, \phi) = \sum_{k=1}^N{c_k\sin(k\alpha\phi)\mathsf{P}^
{-k\alpha}_\nu
(\cos\theta)},\qquad \text{where}\quad\nu(\nu+1)=\lambda
\end{equation}
vanishes numerically on the opposite side of the triangle. This is the
analogue of Equation~\eqref{eq:expansion-u} in the previous section.

\begin{figure}
  \centering
  \includegraphics[trim=50 130 50 130,clip=true,height=6cm]%
  {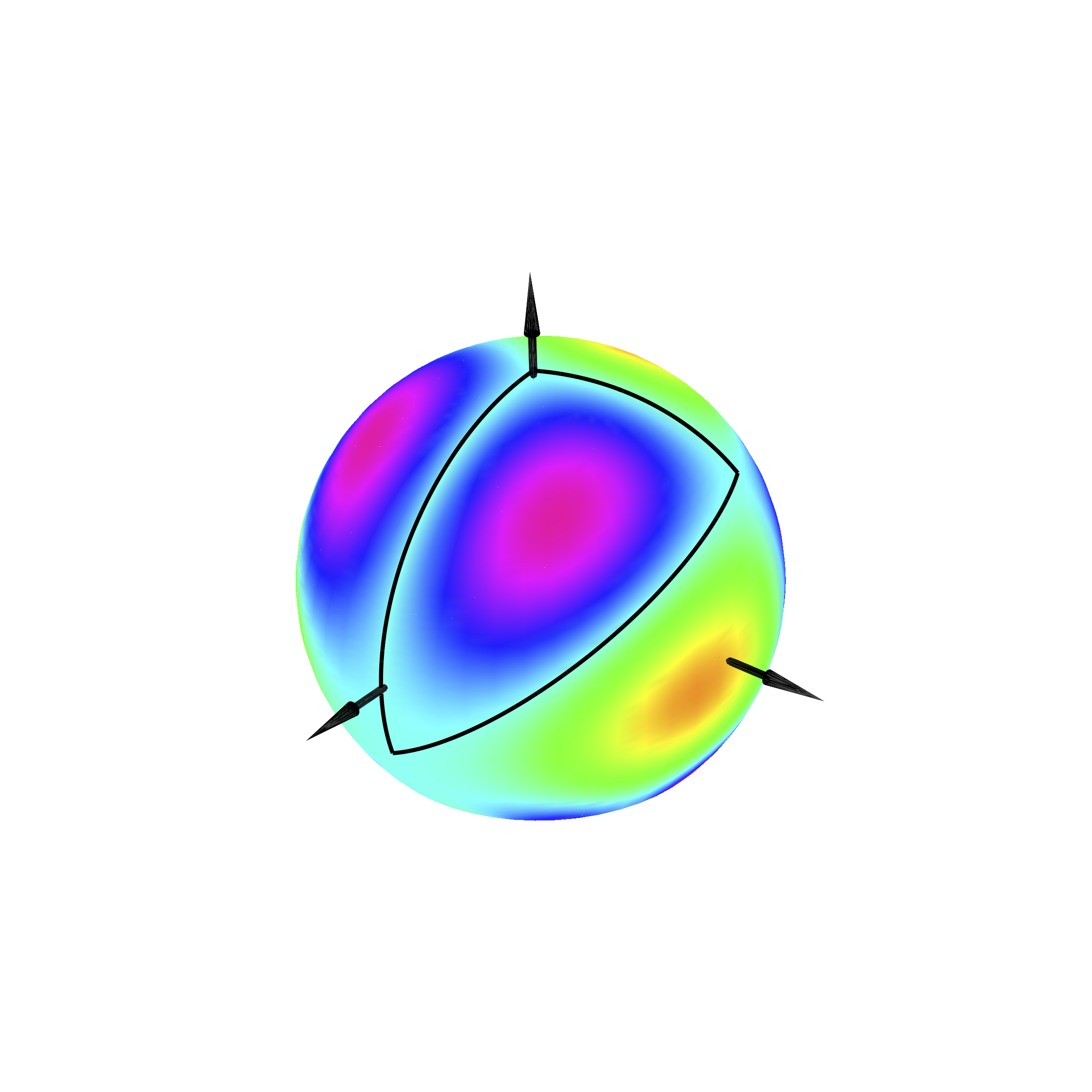}
  \caption{Approximate eigenfunction from Eq.~\eqref{eq:expansion-Ferrers} for the fundamental eigenvalue on
  the spherical triangle $T_2$,
  with angles $2\pi/3$ (at the north pole), $\pi/3$ (on the meridian
  $\phi=0$) and $\pi/2$. The function is defined except at the
  south pole and satisfies~$\Delta
  u+\lambda u=0$ in that domain; it is numerically~0 on the boundary
  of~$T_2$.
  The values of the function are represented by colors. The
  triangle is drawn in black.
  } \label{fig:eigenfunction}
  \end{figure}

\subsection{Upper bound on the maximum}
By design, $u$ is identically zero on two sides of
the triangle and it is sufficient to bound the maximum on the third
side. From the differential equation~\eqref{eq:legendre} one can
obtain a recurrence for the Taylor coefficients of the Ferrers
functions. This allows us to compute Taylor expansions of $u$ when
bounding the maximum. In the computations we use Taylor expansions of
the same order as the number of terms in $u$.

\subsection{Lower bound on the norm}\label{sec:regtriglowerbound}
We compute the norm on the subset of $D$ which in spherical
coordinates is given by the rectangle
\begin{equation*}
  G = \{(\theta, \phi) \in \mathbb{S}^{2}: 0 \leq \phi \leq \pi/\alpha,\ 0 \leq \theta \leq \beta\},
\end{equation*}
where $\beta$ is the minimum value of $\theta$ on the lower
boundary of the triangle. The rectangle $G$ is equal to $D$
precisely when the two angles not at the north pole are both equal to
$\pi/2$. The lower bound is obtained from
\begin{equation*}
\|u\|_2=\int_Du^2\,dx\ge\int_{G}u^{2}\,dx
  = \int_{0}^{\pi/\alpha}\int_{0}^{\beta}u(\theta, \phi)^{2}\sin(\theta)d\theta d\phi.
\end{equation*}
As in the plane, when $u$ is of the form given by
Equation~\eqref{eq:expansion-Ferrers} the orthogonality of the family
$\sin(k\alpha \phi)$ on the interval $[0,\pi/\alpha]$ simplifies this
integral to
\begin{align*}
  \|u\|_2\ge \frac{\pi}{2\alpha}\sum_{k = 1}^{N}c_{k}^{2}
                  \int_{0}^{\beta}\mathsf{P}^{-k\alpha}_{\nu}(\cos
                  \theta)^{2}\sin\theta\,d\theta.
\end{align*}
The integrals
$\int_{0}^{\beta}\mathsf{P}^{-k\alpha}_{\nu}(\cos\theta)^{2}\sin
\theta\,d\theta$ can now be efficiently computed with a certified
integrator, giving us a lower bound for the norm. In practice we
compute the integral from $\epsilon$ to $\beta$ for a small
$\epsilon$ to avoid having to deal with the branch cut at 0, while
not losing too much precision on the bound.

\subsection{Certification of the index}\label{sec:cert_index}
\begin{figure}
  \centering
  \includegraphics[trim=50 90 50 100,clip=true,height=6cm]%
  {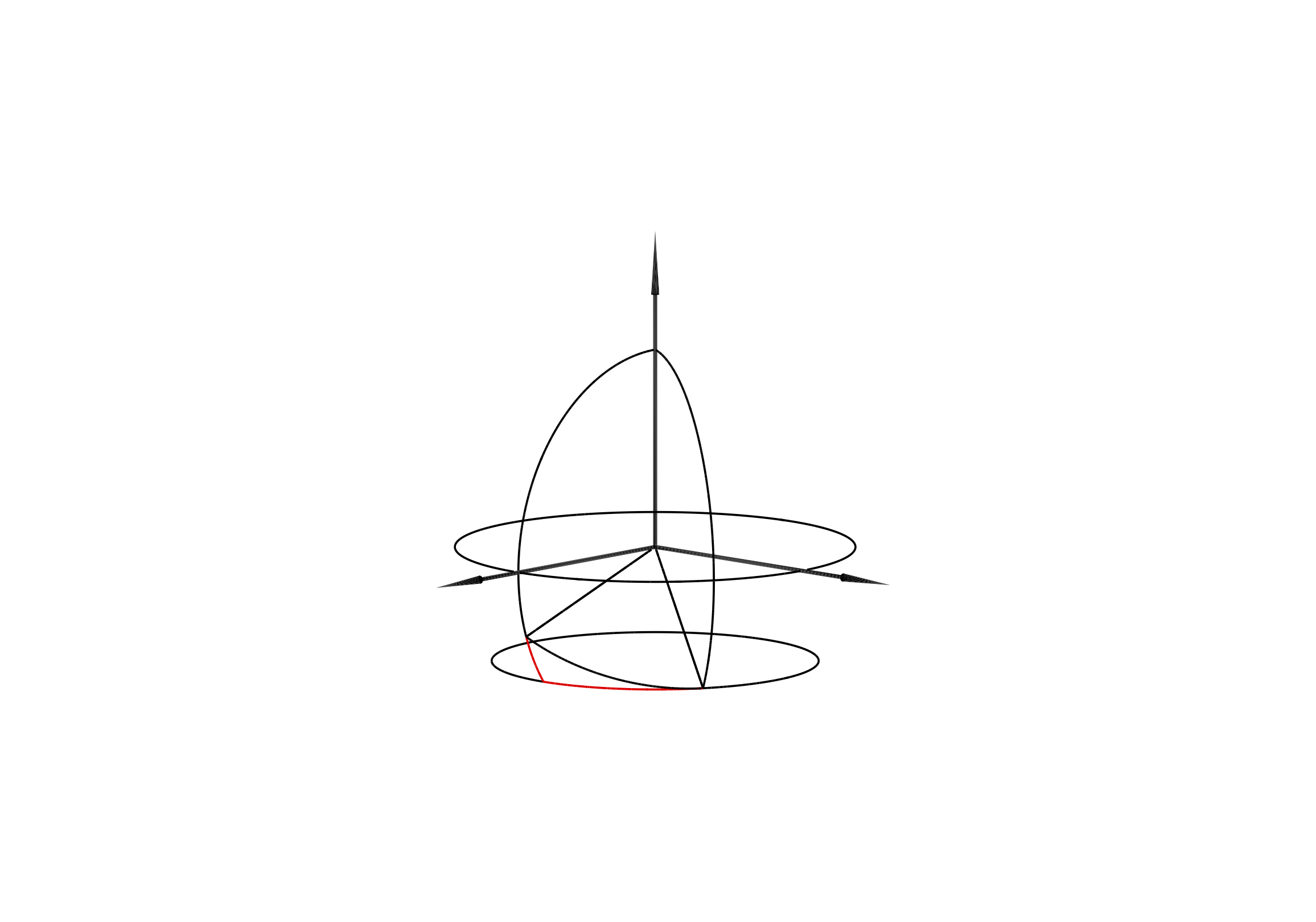}
  \caption{Enclosing sector of a spherical cap (in red) for the
    spherical triangle $T_2$ (in black), with angles $\pi/3$ (at the
    north pole), $2\pi/3$ (on the meridian $\phi=0$) and
    $\pi/2$.} \label{fig:boundtriang}
\end{figure}

That the computed eigenvalue is the fundamental one is certified by
showing that the second eigenvalue of the Laplacian on a larger domain
is larger than it. This implies that this computed eigenvalue is
smaller than the second one and thus has to be the first one.

We consider the domain given by a sector of a spherical cap with its
vertex at the north pole and polar angle $\theta_T$ equal to the
maximal polar angle of the points in the original triangle (see
Figure~\ref{fig:boundtriang}.) The corresponding eigenfunctions are
the products $\sin(k\alpha\phi)\mathsf{P}_\nu^ {-k\alpha}(\cos\theta)$
such that $\mathsf{P}_\nu^ {-k\alpha}(\cos\theta_T)=0$, with
corresponding eigenvalue $\lambda=\nu (\nu+1)$ (see
Equation~\eqref{eq:expansion-Ferrers}). Thus we want to find the
second one in the infinite set of zeros in $\nu$ of the set of Ferrers
functions~$\{ \mathsf{P}_\nu^{-k\alpha}(\cos\theta_T)\mid
k\in\mathbb{N}\setminus \{0\}\}$. The study of this infinite set of
zeroes is simplified by the following.
\begin{lemma}\label{lemma:sturm}
Assume $x_0\in(-1,1)$ is the largest zero of
  $\mathsf{P}^{\mu}_{\nu}(x)$ (with $\mu \leq 0$), then for any
  $(\tilde{\mu},\tilde{\nu})$ such that $\tilde{\mu}<\mu$ and
  $\tilde{\nu}\le\nu$, the function
  $\mathsf{P}^{\tilde{\mu}}_{\tilde{\nu}}(x)$ does not vanish in the
  interval $[x_0,1)$.
\end{lemma}
Thus, denoting by $\zeta_{k,j}$ the $j$th zero of $\mathsf{P}^
{-k\alpha}_\nu
(\cos\theta_T)$ as a function of~$\nu$, if a triangle has computed
eigenvalue
$\lambda = \nu(\nu + 1)$, it is sufficient to show that
$\zeta_{1, 2}$ and $\zeta_{2, 1}$ are both larger than $\nu$ to
certify that $\lambda$ is indeed the fundamental eigenvalue.
\begin{proof}[Proof of the lemma]The basic idea is to use Sturm's
  comparison theorem in order to show that the largest zero
  $\tilde{x}_0$ of $\mathsf{P}^ {\tilde{\mu}}_{\tilde{\nu}}(x)$ in
  $(-1,1)$ satisfies $ \tilde{x}_0<x_0$.

  The associated Legendre equation~\eqref{eq:legendre} can be
  rewritten
  \[((1-x^2)w')'+q_{\mu,\nu}(x)w=0,\qquad\text{with}\quad q_{\mu,\nu}
    (x)=\nu(\nu+1)-\frac{\mu^{2}}{1-x^{2}}\] and the inequality
  $q_{\tilde{\mu},\tilde{\nu}}(x)<q_{\mu,\nu}(x)$ for $x\in(-1,1)$
  follows from the hypotheses. The proof is by contradiction.
  Assume~$\mathsf{P}^{\mu}_{\nu}$ does not vanish in
  $ [\tilde{x}_0,1)$. For simplicity of notation, write
  $w=\mathsf{P}^ {\mu}_{\nu}$ and
  $\tilde{w}=\mathsf{P}^ {\tilde{\mu}}_{\tilde{\nu}}$. Then a direct
  verification shows Picone's identity
  \[
    \left((1-x^2)\frac{\tilde{w}}w(
      \tilde{w}'w-w'\tilde{w})\right)'=\\
    (q_{\mu,\nu}-q_{\tilde{\mu},\tilde{\nu}})
    \tilde{w}^2+(1-x^2)\left(\tilde{w}'-w'\frac{\tilde{w}} {w}\right)^2,
  \]
  whose right-hand side is positive in $[\tilde{x}_0,1)$. This
  implies that the function
  \[(1-x^2)\frac{\tilde{w}}w(\tilde{w}'w-w'\tilde{w})\] is increasing
  in that interval. However, it is 0 at~$\tilde{x}_0$ while as
  $x\rightarrow 1-$, it behaves like
  \[(\tilde{\mu}-\mu)2^{\tilde{\mu}}(1-x)^{-\tilde{\mu}}\] and thus
  tends to~0 at~$1$, a contradiction.
\end{proof}
Using standard interval methods we can isolate all the roots of the
first and second Ferrers functions in the interval $[0, \nu]$ and if
there is exactly one root then we are sure that the second ones,
$\zeta_ {1, 2}$ and $\zeta_{2, 1}$ are larger than $\nu$ and therefore
that $\lambda$ is the fundamental eigenvalue.

The method can fail in case either of $\zeta_{1, 2}$ and
$\zeta_{2, 1}$ is less than $\nu$. In that case, we cannot
conclude anything about the index of $\lambda$. One option then is
to try a different orientation of the triangle. Depending on which
angle of the spherical triangle is placed at the north pole, we get a
different set of Ferrers functions to consider. In several of the
examples below the method fails with the original orientation but
there is always at least one orientation in which it succeeds.

The same approach can be used to find the spherical cap with
fundamental eigenvalue $\lambda$, required for the Faber-Krahn
inequality in Equation \eqref{eq:faber-krahn}. For a spherical cap,
the eigenfunctions are $\sin(m\phi)\mathsf{P}_\nu^{m}(\cos\theta)$,
with $m \in \mathbb{Z}$ for it to be continuous on the whole cap, and
with $\mathsf{P}_\nu^{m}(\cos\theta) = 0$ along the boundary of the
cap. The corresponding eigenvalue is $\nu(\nu + 1)$. The polar angle
of the cap with fundamental eigenvalue $\lambda$ is thus determined by
finding the largest zero $x_{0} \in (-1, 1)$ of
$\mathsf{P}^{m}_{\nu}(x)$ with $\lambda = \nu(\nu + 1)$ and
$m \in \mathbb{Z}$. Lemma~\ref{lemma:sturm} reduces the computation
to the case \(m = 0\), a Legendre function. The area of the cap is
then given by $2\pi(1 -
x_{0})$
and during the computation of the norm, as described in
Section~\ref{section:norm}, any region of area larger than this will
need to be split.

\subsection{Regular triangles}\label{sec:regular}
We now present the results obtained by this method for the triangles
appearing in Table 3 in the work of Bogosel \emph{et al.}
\cite{BogoselPerrollazRaschelTrotignon2020} which have at most one
singular vertex, i.e., a vertex whose angle is not of the form
$\pi/k$ for some integer $k$, and for which the eigenvalue is not
exactly known. These triangles are given in
Table~\ref{table:regular-triangles}, together with their computed
eigenvalues.

We start by showing the successive steps of the method for the
triangle~$T_{2}$ with angles~$(2\pi/3,\pi/3,\pi/2)$ (see
Figure~\ref{fig:eigenfunction}).

\begin{table}
  \centering
  \begin{tabular}{|c|c|c|c|c|c|}
    \hline
     &  Number in \cite{BogoselPerrollazRaschelTrotignon2020} & Angles
     & Eigenvalue & Eigenvalue in
     \cite{BogoselPerrollazRaschelTrotignon2020} \\ \hline
    $T_{1}$ & 8  & $\left(\frac{3\pi}{4}, \frac{\pi}{3}, \frac{\pi}{2}\right)$
                                             & 12.400051652843377905
                                                          & 12.400051 \\ \hline
    $T_{2}$ & 9  & $\left(\frac{2\pi}{3}, \frac{\pi}{3}, \frac{\pi}{2}\right)$
                                             & 13.744355213213231835
                                                          & 13.744355 \\ \hline
    $T_{3}$ & 11 & $\left(\frac{2\pi}{3}, \frac{\pi}{4}, \frac{\pi}{2}\right)$
                                             & 20.571973537984730557
                                                          & 20.571973 \\ \hline
    $T_{4}$ & 12 & $\left(\frac{2\pi}{3}, \frac{\pi}{3}, \frac{\pi}{3}\right)$
                                             & 21.309407630190445259
                                                          & 21.309407 \\ \hline
    $T_{5}$ & 13 & $\left(\frac{3\pi}{4}, \frac{\pi}{4}, \frac{\pi}{3}\right)$
                                             & 24.456913796299111694
                                                          & 24.456913 \\ \hline
    $T_{6}$ & 16 & $\left(\frac{2\pi}{3}, \frac{\pi}{4}, \frac{\pi}{4}\right)$
                                             & 49.109945263284609920
                                                          & 49.109945 \\ \hline
  \end{tabular}
  \caption{Triangles with at most one singular vertex from Table 3 in
    \cite{BogoselPerrollazRaschelTrotignon2020}. Certified, correctly
    rounded, 20 digit eigenvalues are given. Previously computed
    eigenvalues are also shown. See
    Table~\ref{table:triangles-high-precision-results} for more
    digits.}
  \label{table:regular-triangles}
\end{table}

\begin{figure}
  \centering
  \begin{subfigure}[t]{0.45\textwidth}
    \includegraphics[width=\textwidth]{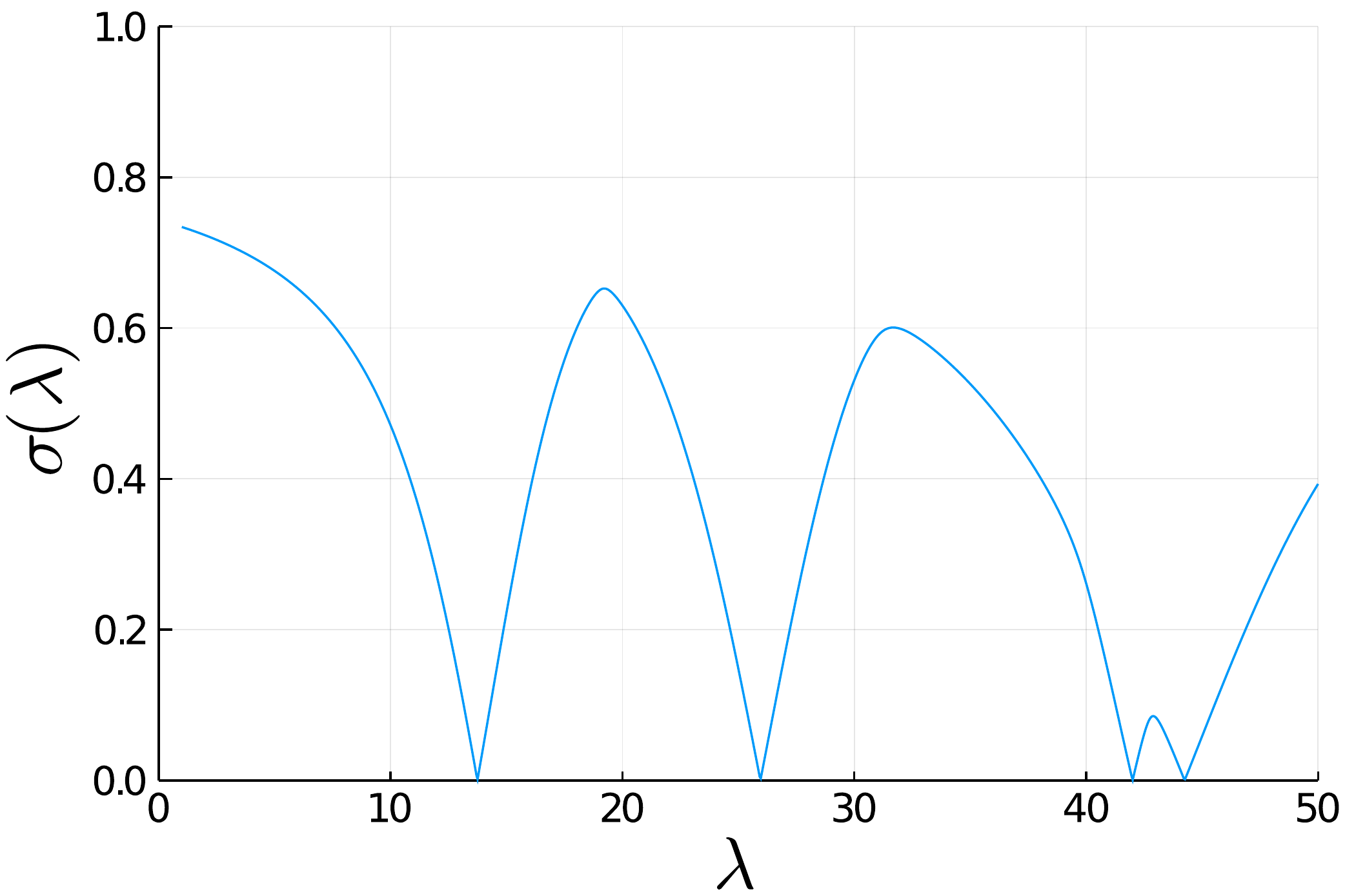}
    \caption{The function $\sigma(\lambda)$.}
    \label{fig:regular-triangles-example-sigma}
  \end{subfigure}
  \hspace{0.05\textwidth}
  \begin{subfigure}[t]{0.45\textwidth}
    \includegraphics[width=\textwidth]{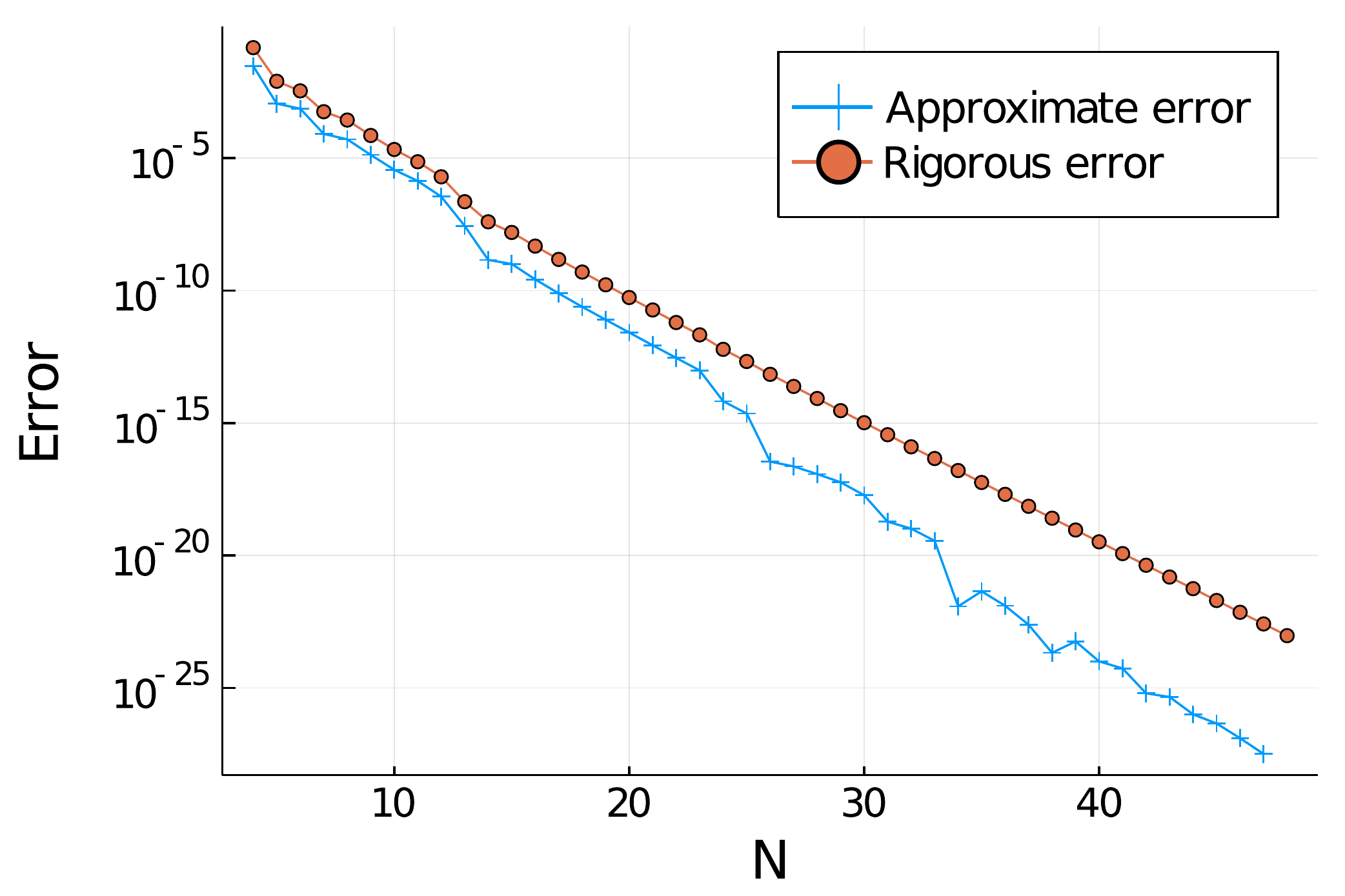}
    \caption{Convergence to the eigenvalue corresponding to the first
      minimum of $\sigma(\lambda)$.}
    \label{fig:regular-triangles-example-convergence}
  \end{subfigure}
  \caption{Results for the triangle $T_{2}$ from
    Table~\ref{table:regular-triangles}}
\end{figure}

The value of $\sigma(\lambda)$ is given in
Figure~\ref{fig:regular-triangles-example-sigma}. It was generated
using an expansion at the singular vertex, with 8 terms,
using 16 random points in the interior and 16 points on the boundary
opposite to the singular vertex. The plot shows 4 minima, the first of
which we need to certify as corresponding to the fundamental
eigenvalue.

The convergence towards the first minimum is shown in
Figure~\ref{fig:regular-triangles-example-convergence}. Similarly to
Figure~\ref{fig:lshape-convergence}, it shows two kinds of
convergence, the approximate error computed as the difference to the
``exact'' solution obtained with $N = 48$ and the certified
enclosure. The quotient between the approximate and the certified
error increases slightly with $N$, the precision ``lost'' when going
from the approximate error to the certified error is slightly more
than 5 digits. For $N = 48$ we get the enclosure
$\lambda \in [13.7443552132132318354011 \pm 3.11\,10^{-23}]$. See
Table~\ref{table:triangles-high-precision-results} for more digits.

With the vertex with angle $2\pi/3$ placed at the north pole we get
the zeros $\zeta_{1, 2} \in [3.6550969 \pm 4.82\,10^{-8}]$ and
$\zeta_{2, 1} \in [3.4315893 \pm 5.43,10^{-8}]$ for the enclosing
spherical cap sector. Since
$\lambda < \zeta_{1, 2}(\zeta_{1, 2} + 1)$ and
$\lambda < \zeta_{2, 1}(\zeta_{2, 1} + 1)$ it is smaller than the
second eigenvalue of the enclosing spherical cap sector and must
therefore correspond to the fundamental eigenvalue of the triangle.

\begin{figure}
  \centering
  \includegraphics[height=5cm]{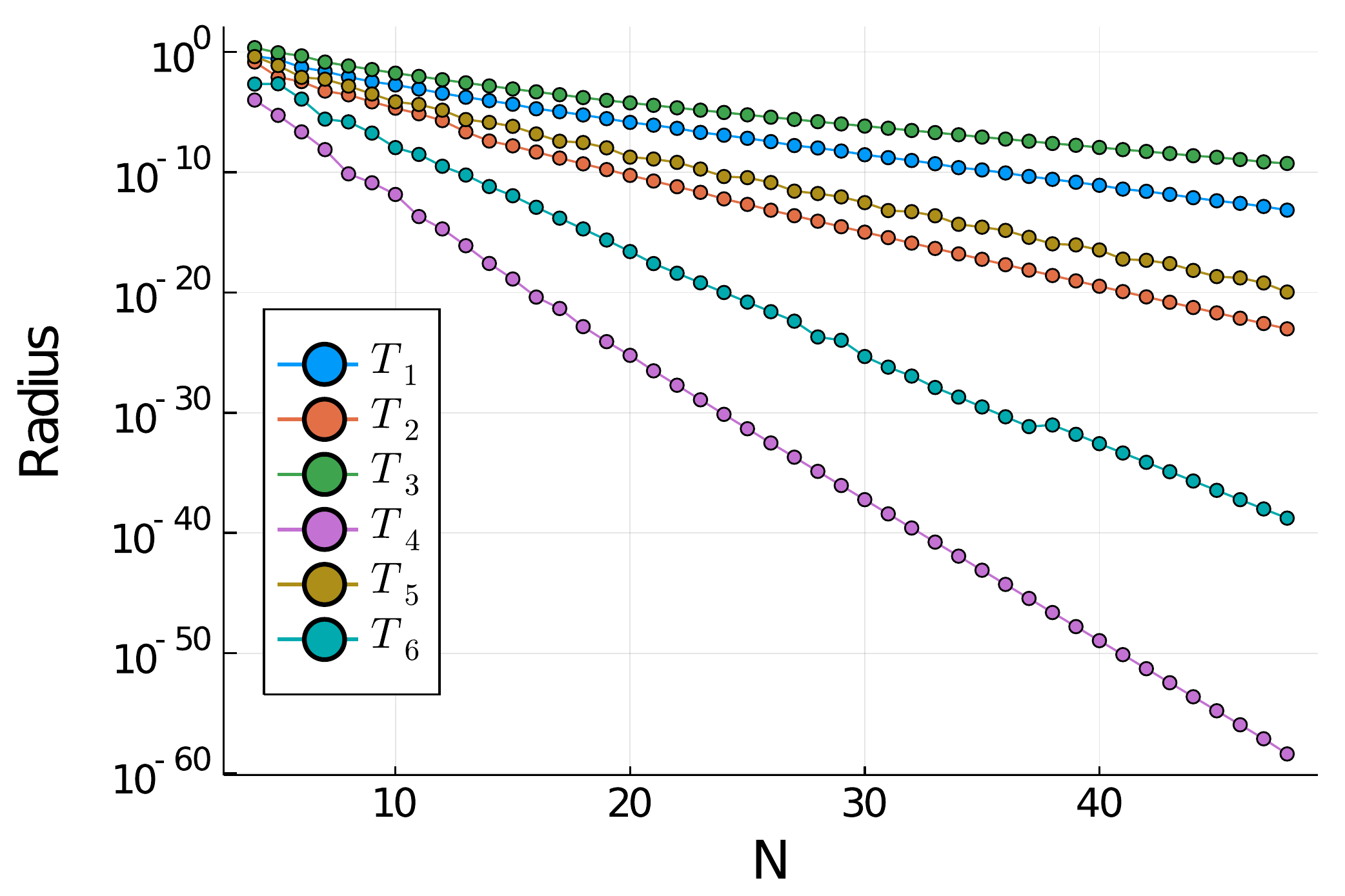}
  \caption{Convergence of the enclosure for the six triangles in
    Table~\ref{table:regular-triangles}.}
  \label{fig:regular-triangles-enclosure}
\end{figure}

For the other triangles in Table~\ref{table:regular-triangles}, the
method is the same except for the triangles $T_{4}$ and $T_{6}$.
For these triangles the two non-singular angles are the same. This
symmetry in the domain implies a corresponding symmetry for the
eigenfunction. The approximate solution can be forced to have the same
symmetry by using only every second term from the sum in
Equation~\eqref{eq:expansion-Ferrers}, which improves the convergence
rate. Figure~\ref{fig:regular-triangles-enclosure} shows the
convergence for the radius of the enclosures. The rate of convergence
varies between the triangles, the best convergence being obtained for
the triangles $T_{4}$ and $T_{6}$ where the mentioned symmetry was
used. Even though they converge at different rates they all show
linear convergence. As in the case of the L-shaped domain the
condition
number of $R(\lambda)$ is huge for all of them, $T_{6}$ has the
highest value at $10^{300}$. As in the previous cases the solutions is
however still small on the boundary.

For all these triangles we are also able to certify that the computed
eigenvalue indeed corresponds to the fundamental eigenvalue by lower
bounding the second eigenvalue of the enclosing circular cap sector.

\subsubsection*{High precision results}
\label{sec:high-prec-results}
\begin{figure}
  \centering
  \begin{subfigure}[t]{0.48\textwidth}
    \includegraphics[width=\textwidth]{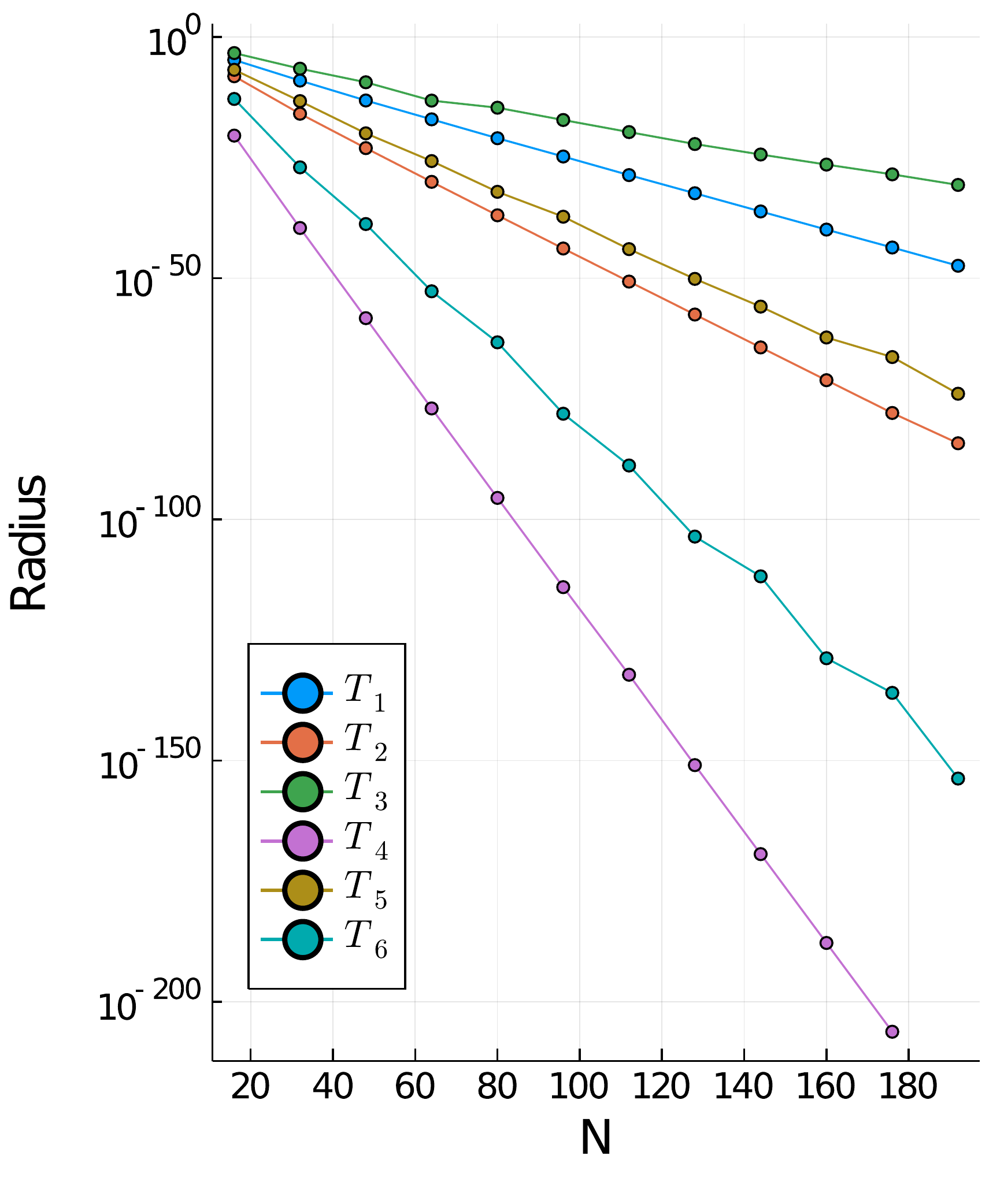}
    \caption{High precision computations of the enclosure}
    \label{fig:regular-triangles-high-precision-convergence}
  \end{subfigure}
  \hspace{0.02\textwidth}
  \begin{subfigure}[t]{0.48\textwidth}
    \includegraphics[width=\textwidth]{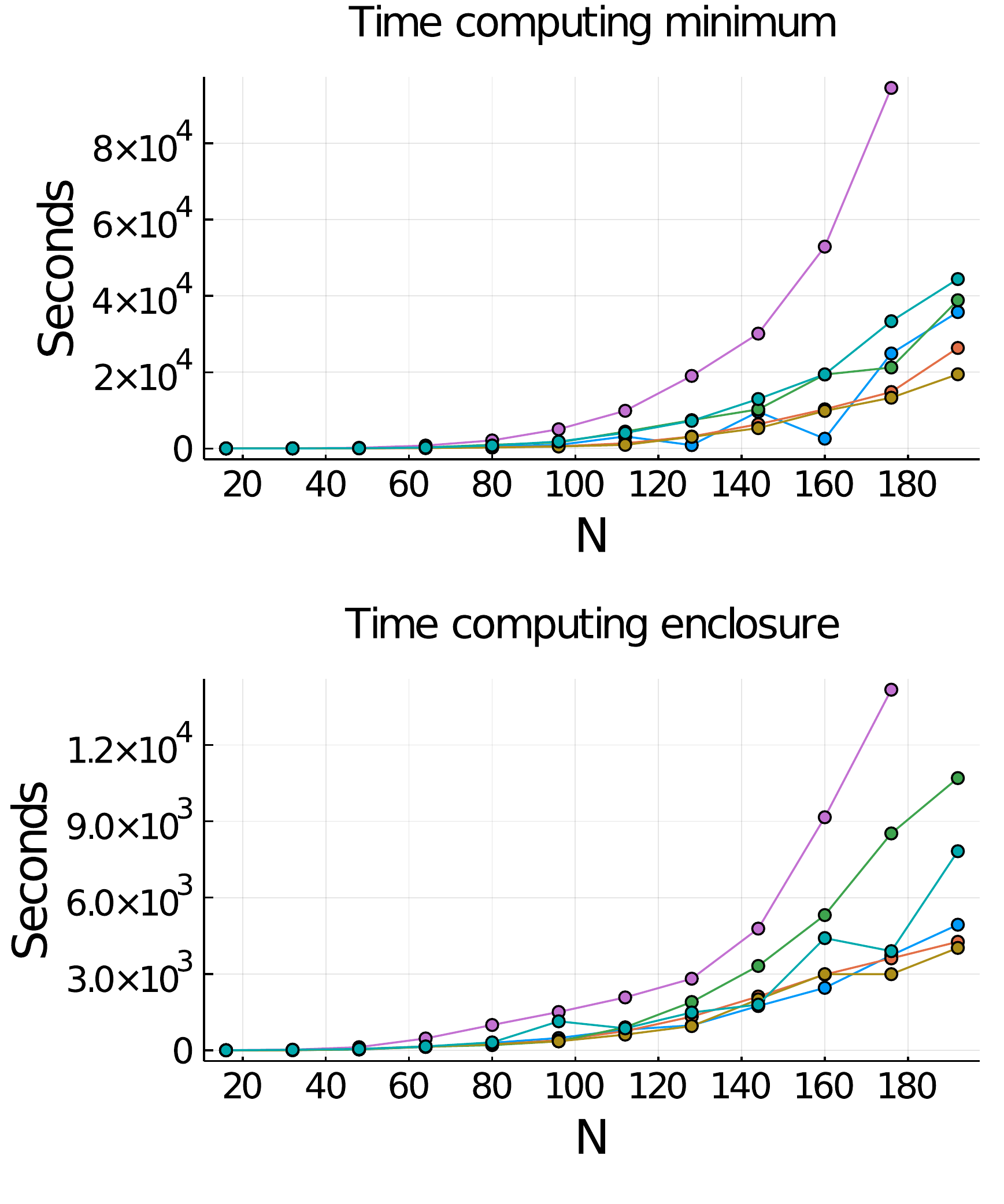}
    \caption{Time for computing the minimum of $\sigma(\lambda)$ as
      well as the enclosure.}
    \label{fig:regular-triangles-high-precision-time}
  \end{subfigure}
  \label{fig:regular-triangles-high-precision}
  \caption{High precision computations for the triangles in
    Table~\ref{table:regular-triangles}.}
\end{figure}

We give high precision computations of the same eigenvalues in
Figure~\ref{fig:regular-triangles-high-precision-convergence}.
Compared to Figure~\ref{fig:regular-triangles-enclosure} we start the
computations at $N = 16$ and increase $N$ by 16 at a time.
Figure~\ref{fig:regular-triangles-high-precision-time} shows the time
taken for both computing the minimum of $\sigma(\lambda)$ and
computing the enclosure for different triangles as $N$ varies. This
shows that for large values of $N$ the cost of certifying the
enclosure is a relatively small part of the computation, and in
particular, the time for bounding the norm accounts for only a few
seconds of the total time. Correctly rounded eigenvalues are given in
Table~\ref{table:triangles-high-precision-results}.

\subsection{Singular triangles}\label{sec:singular-triangles}
As shown above, the method of particular solutions works well for the
triangles in Table~\ref{table:regular-triangles}, using an expansion
at the single singular vertex. This is not sufficient for the
singular triangles listed in Table~\ref{table:singular-triangles}.
Those are the triangles from Table~3 in~%
\cite{BogoselPerrollazRaschelTrotignon2020} with more than one
singular vertex.

\begin{table}
  \centering
  \begin{tabular}{|c|c|c|c|c|}
    \hline
     & Number in \cite{BogoselPerrollazRaschelTrotignon2020} & Angles
     & Eigenvalue & Eigenvalue in
     \cite{BogoselPerrollazRaschelTrotignon2020} \\ \hline
    $T_{7}$ & 1  & $\left(\frac{2\pi}{3}, \frac{3\pi}{4}, \frac{3\pi}{4}\right)$
                                             & 4.2617347552939870857
                                                          & 4.261734 \\ \hline
    $T_{8}$ & 2  & $\left(\frac{2\pi}{3}, \frac{2\pi}{3}, \frac{2\pi}{3}\right)$
                                             & 5.1591456424665417112
                                                          & 5.159145 \\ \hline
    $T_{9}$ & 3 & $\left(\frac{\pi}{2}, \frac{2\pi}{3}, \frac{3\pi}{4}\right)$
                                             & 6.2417483307263342368
                                                          & 6.241748 \\ \hline
    $T_{10}$ & 4 & $\left(\frac{\pi}{2}, \frac{2\pi}{3}, \frac{2\pi}{3}\right)$
                                             & 6.7771080545983009574
                                                          & 6.777108 \\ \hline
  \end{tabular}
  \caption{Triangles with more than one singular vertex from Table 3
    in \cite{BogoselPerrollazRaschelTrotignon2020}. Certified,
    correctly rounded, 20 digit eigenvalues are given. See
    Table~\ref{table:triangles-high-precision-results} for more
    digits.}
  \label{table:singular-triangles}
\end{table}

\begin{figure}
  \centering
  \begin{subfigure}[t]{0.45\textwidth}
    \includegraphics[width=\textwidth]{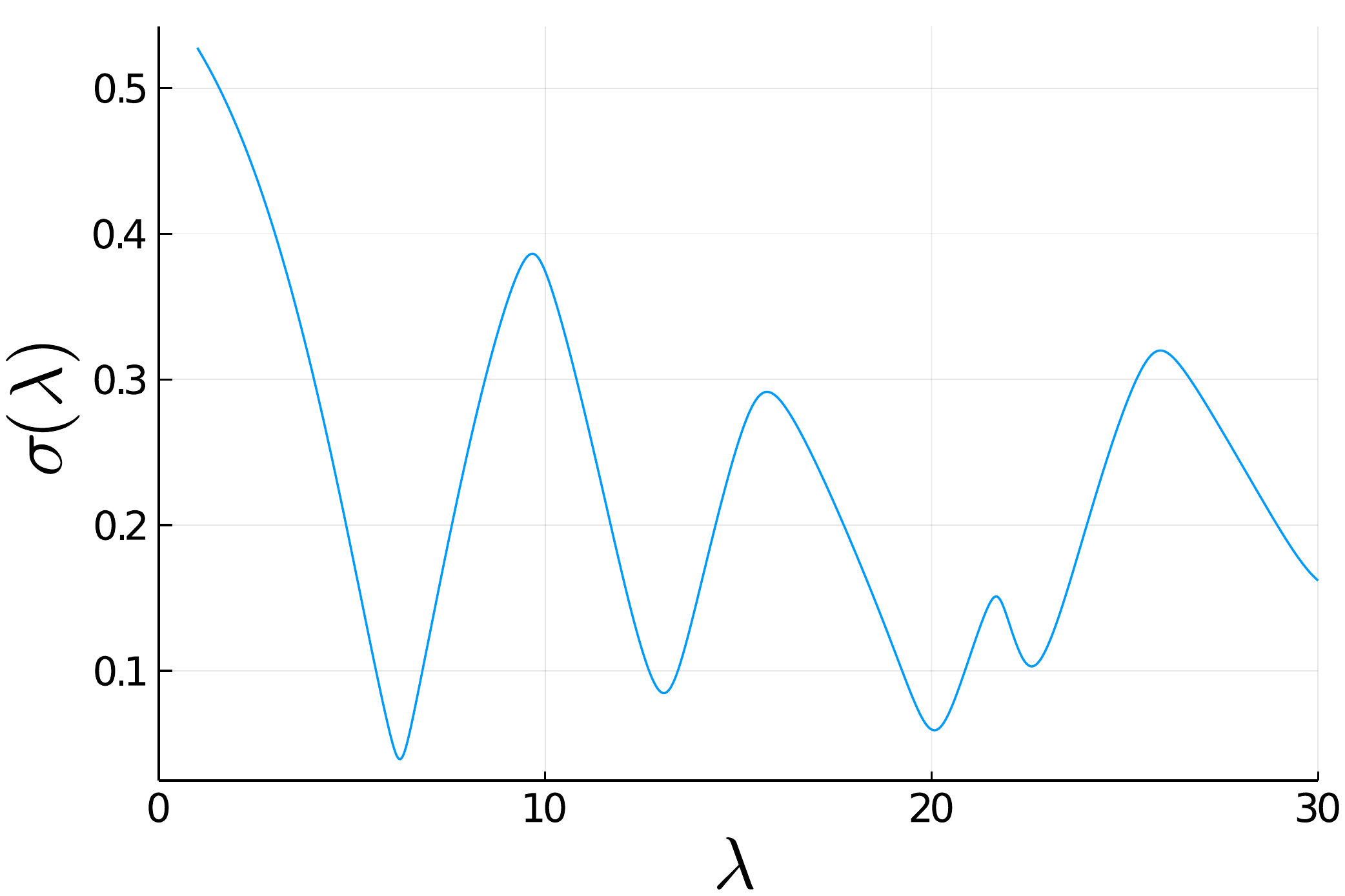}
    \caption{The function $\sigma(\lambda)$ with $N = 16$ terms.}
    \label{fig:singular-triangles-example-sigma-withoutinterior}
  \end{subfigure}
  \hspace{0.05\textwidth}
  \begin{subfigure}[t]{0.45\textwidth}
    \includegraphics[width=\textwidth]{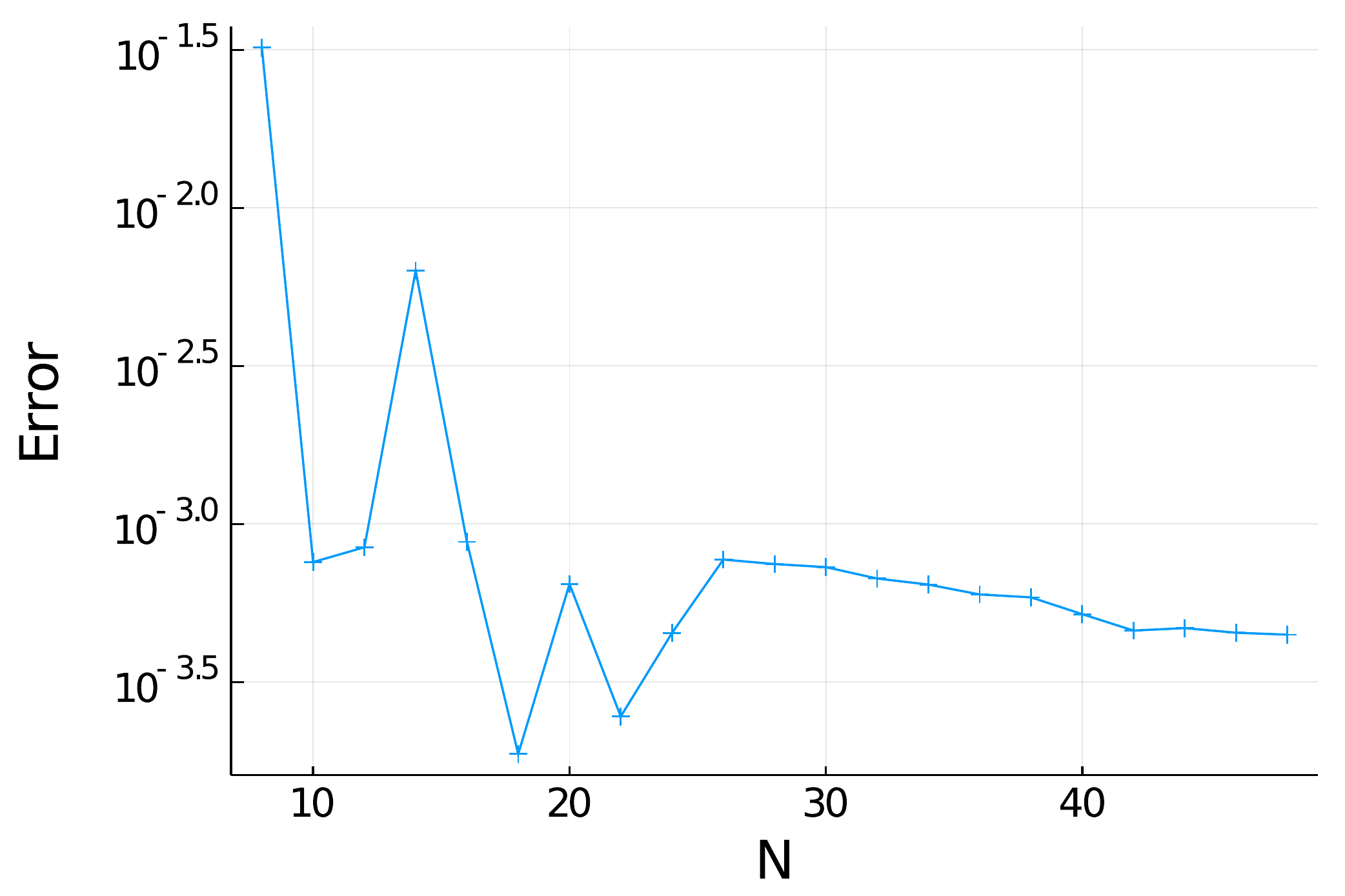}
    \caption{Convergence of the first minimum. No improvement is seen
      as the number of terms is increased past 26.}
    \label{fig:singular-triangles-example-approximate-convergence-withoutinterior}
  \end{subfigure}
  \caption{Results for the triangle $T_{9}$ from
    Table~\ref{table:singular-triangles} using expansions with
    $N/2$ terms at each of the two singular vertices.}
\end{figure}

For such singular cases, a solution suggested by Betcke and Trefethen
is to use expansions at all singular vertices. For a spherical
triangle where all of the vertices are singular our candidate
eigenfunction will be of the form
\begin{equation*}
  u(\theta, \phi) = u_{1}(\theta_{1}, \phi_{1}) + u_{2}(\theta_{2},
  \phi_{2}) + u_{3}(\theta_{3}, \phi_{3})
\end{equation*}
with
\[u_{i}(\theta, \phi) = \sum_{k = 1}^{N_{i}} c_{i,k}
  \sin(k\alpha_{i}\phi) \mathsf{P}^{-k\alpha_{i}}_\nu(\cos\theta)\]
and $(\theta_{i}, \phi_{i})$ corresponds to $(\theta, \phi)$ in
spherical coordinates with vertex $i$ on the north pole. This method
does not work well in our cases. Triangle $T_{9}$ from
Table~\ref{table:singular-triangles} has two singular vertices. Using
expansions with 8 terms at each one gives the plot of
$\sigma(\lambda)$ seen in
Figure~\ref{fig:singular-triangles-example-sigma-withoutinterior}.
This plot is less smooth than the corresponding one for the regular
triangle (Figure~\ref{fig:regular-triangles-example-sigma}). As the
number of terms is increased the picture does not improve. Convergence
of the first minimum is show in
Figure~\ref{fig:singular-triangles-example-approximate-convergence-withoutinterior}
where there error is computed by comparing it to
$\lambda = 6.24174833072633424$. Increasing the number of terms past
26 does not improve the approximation.

\begin{figure}
  \centering
  \begin{subfigure}[t]{0.45\textwidth}
    \includegraphics[width=\textwidth]{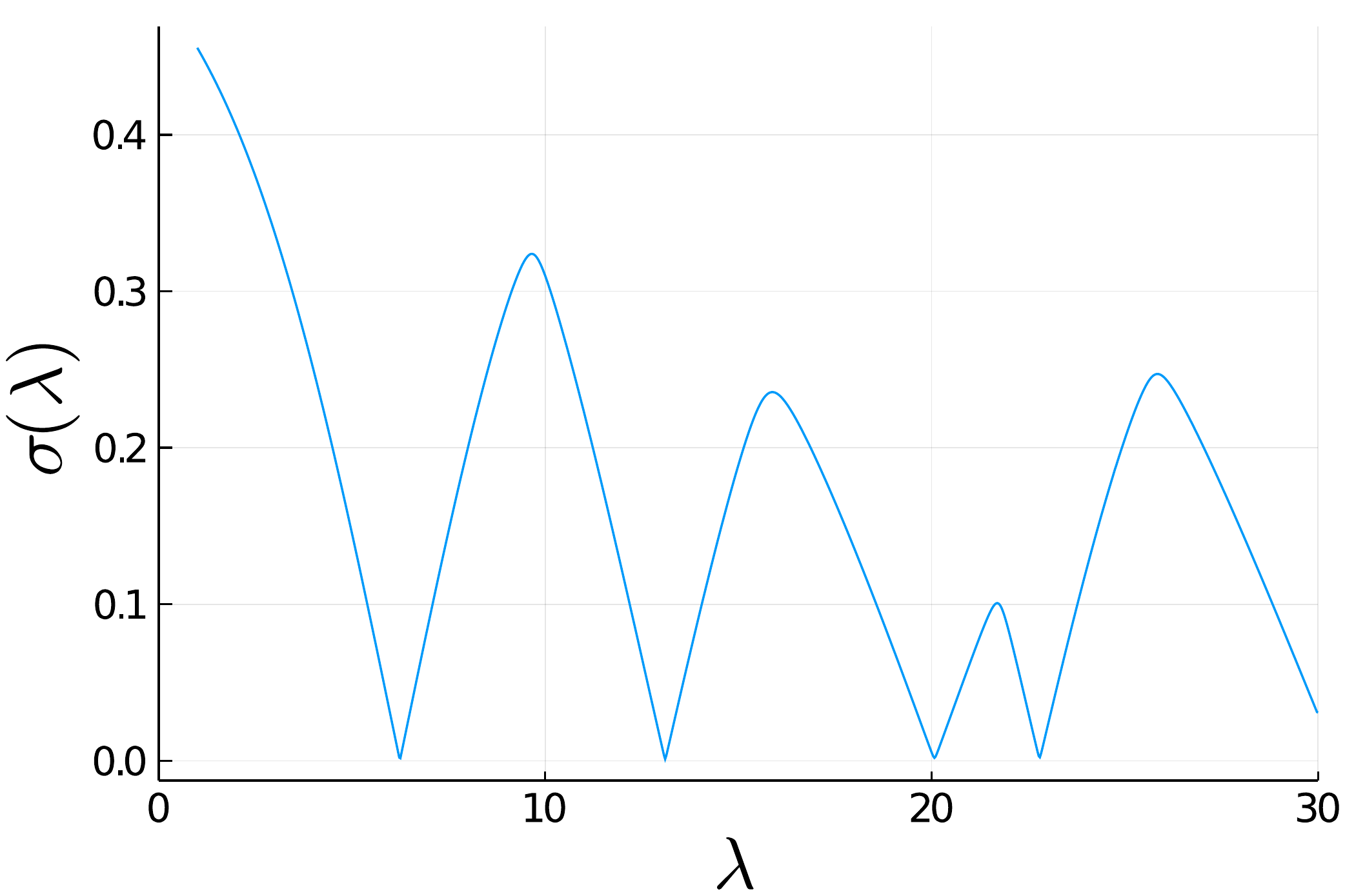}
    \caption{The function $\sigma(\lambda)$ using expansions with
      3 terms at each vertex as well as 12 interior terms.}
    \label{fig:singular-triangles-example-sigma-withinterior}
  \end{subfigure}
  \hspace{0.05\textwidth}
  \begin{subfigure}[t]{0.45\textwidth}
    \includegraphics[width=\textwidth]{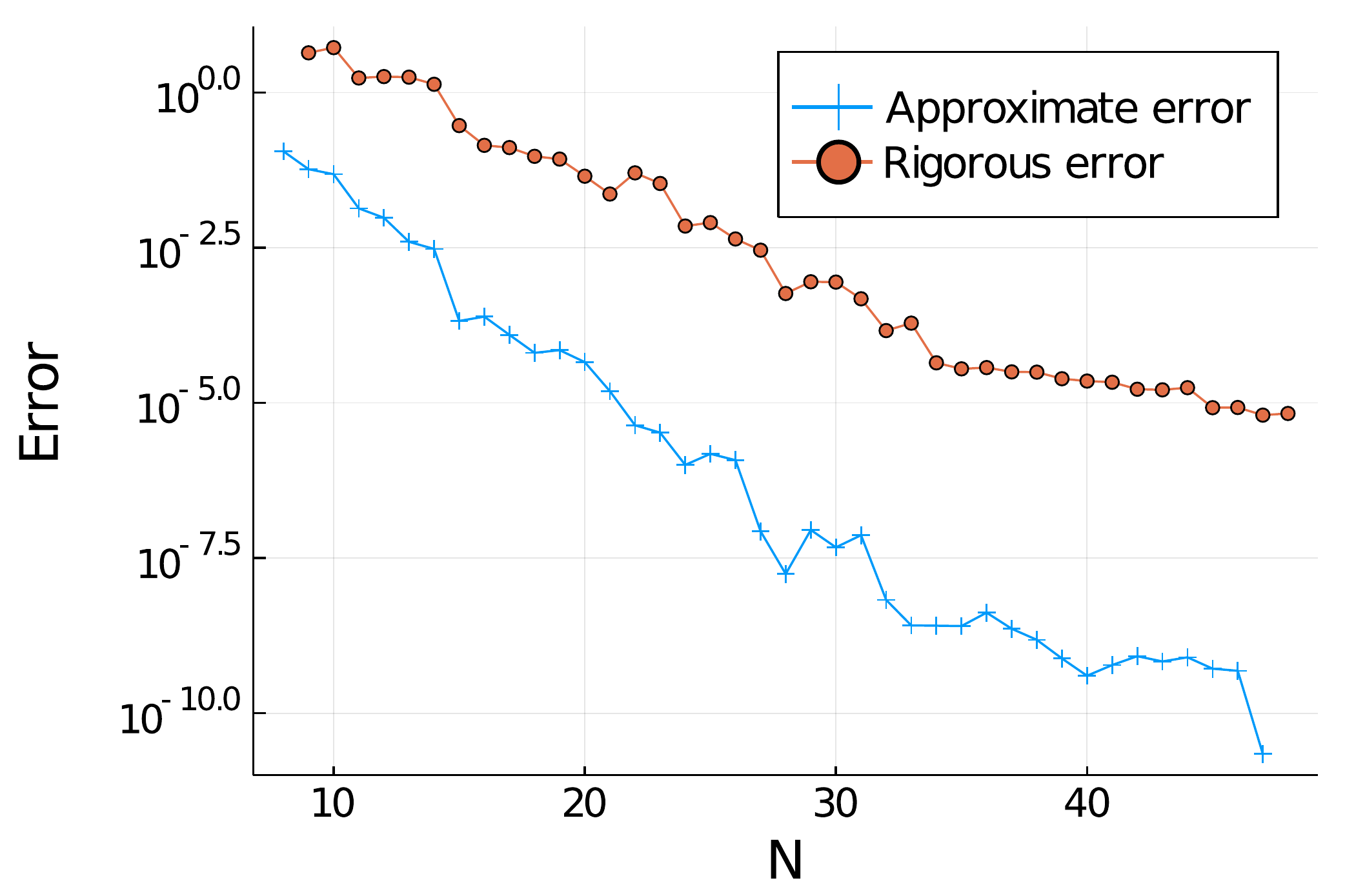}
    \caption{Convergence towards the eigenvalue corresponding to the
      first minimum for triangle $T_{9}$ from
      Table~\ref{table:singular-triangles}.}
    \label{fig:singular-triangles-example-convergence-withinterior}
  \end{subfigure}
  \caption{Results for the triangle $T_{9}$ from
    Table~\ref{table:singular-triangles}. .}
\end{figure}

A solution we found to work satisfactorily is to complement the above
approach with an expansion at an interior point. For a triangle with
three singular vertices this gives us the candidate
\begin{equation*}
  u(\theta, \phi) = u_{1}(\theta_{1}, \phi_{1}) + u_{2}(\theta_{2},
  \phi_{2}) + u_{3}(\theta_{3}, \phi_{3})
  + u_{\text{int}}(\theta_{\text{int}}, \phi_{\text{int}})
\end{equation*}
where $(\theta_{\text{int}}, \phi_{\text{int}})$ is given in
spherical coordinates with the interior point placed at the north pole
and $u_{\text{int}}$ contains $N_{\text{int}}$ terms and is of the
form
\begin{align*}
  u_{\text{int}}(\theta, \phi) =& c_{\text{int},1}\mathsf{P}_{\nu}^{0}(\cos(\theta))
  + c_{\text{int},2}\sin(\phi)\mathsf{P}_{\nu}^{1}(\cos(\theta))
  + c_{\text{int},3}\cos(\phi)\mathsf{P}_{\nu}^{1}(\cos(\theta))\\
  &+ c_{\text{int},4}\sin(2\phi)\mathsf{P}_{\nu}^{2}(\cos(\theta))
  + c_{\text{int},5}\cos(2\phi)\mathsf{P}_{\nu}^{2}(\cos(\theta))
  + \cdots.
\end{align*}
Using this expansion for triangle $T_{9}$ with the interior point
chosen to be the center of the triangle, given by the sum of the
vertices normalized to be on the sphere, and 3 terms from the two
singular vertices combined with 12 terms from the interior gives the
plot of $\sigma(\lambda)$ seen in
Figure~\ref{fig:singular-triangles-example-sigma-withinterior}.
Compared to
Figure~\ref{fig:singular-triangles-example-sigma-withoutinterior} the
minima are more distinct. The convergence towards the first minimum is
shown in
Figure~\ref{fig:singular-triangles-example-convergence-withinterior}.
Like Figure~\ref{fig:regular-triangles-example-convergence} it shows
both the approximate error and the radius of the computed enclosure.
Even though there are more oscillations than for the regular
triangles, we see convergence as $N$ is increased. For $N = 48$ we get
the enclosure $\lambda \in [6.24175 \pm 8.42\,10^{-6}]$. We have used
the same number of terms for each of the two singular vertices and
four times that number of terms for the interior point. The upper
bound of the maximum on the boundary is computed in the same way as
for the regular triangles except that it is no longer identically
equal to zero on two of the boundaries. For the lower bound of the
norm, the terms in the expansion are no longer orthogonal and we have
to resort to the second method described in
Section~\ref{section:norm}. Since the first eigenfunction has constant
sign and is very smooth we do not need a very fine partitioning to get
a good lower bound. The interior domain $\Omega'$ we use is the
triangle with vertices given by the points in the middle between the
center and the vertices of the original triangle. Partitioning
$\Omega'$ into four separate triangles yields a sufficiently good
lower bound. Still the computation of the norm is much more costly
than in the case of regular triangles.

\begin{figure}
  \centering
  \includegraphics[width=.6\textwidth]
  {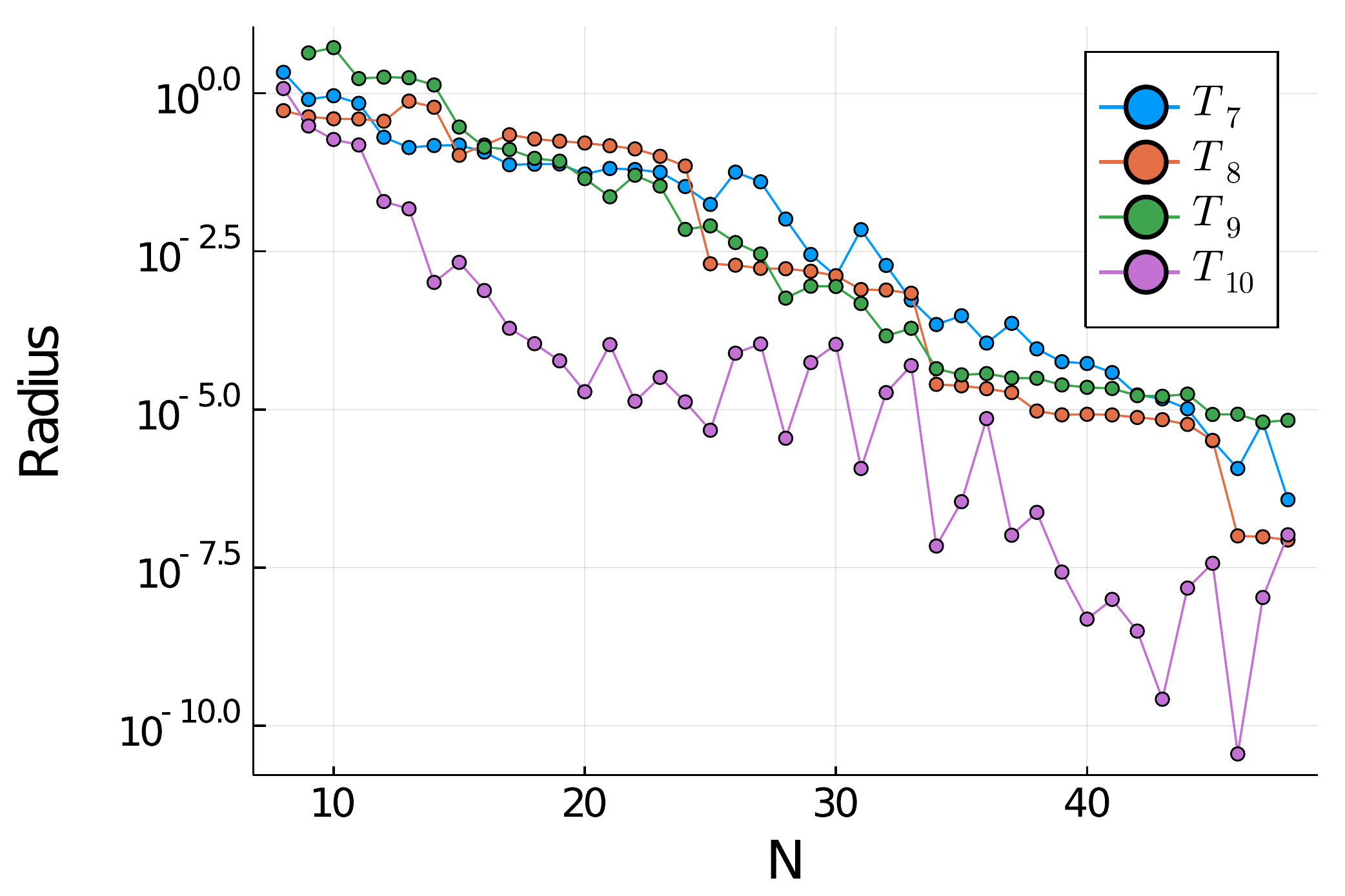}
  \caption{Convergence of the enclosure for the four triangles in
    Table~\ref{table:singular-triangles}.}
  \label{fig:singular-triangles-enclosure}
\end{figure}

Figure~\ref{fig:singular-triangles-enclosure} shows the convergence
for all the triangles in Table~\ref{table:singular-triangles}. Again
there are more oscillations than for the regular triangles but still
convergence as the number of terms is increased. The triangles
$T_{7}$, $T_{8}$ and $T_{10}$ all have symmetries that have
been used to improve the convergence. For $T_{7}$ only every second
term occurs in the expansion for both $u_{1}$ and
$u_{\text{int}}$. In $T_{8}$ we take only every second term from
$u_{1}$, $u_{2}$ and $u_{3}$. Finally for $T_{10}$ every
second term is used in $u_{\text{int}}$. In several cases, there is
potential to make more use of symmetries; for $T_{8}$, this is done
in Section~\ref{sec:3dkreweras}. As in the previous examples
$A(\lambda)$ is ill-conditioned, though slightly less so than before,
with condition number varying between $10^{19}$ and $10^{51}$ for the
four triangles. Again, the solution is still small on the boundary and
we get good enclosures.

As in the case of regular triangles, we are also able to certify that
the computed eigenvalue indeed corresponds to the fundamental
eigenvalue by lower bounding the eigenvalue of the enclosing circular
cap sector.

\begin{figure}
  \centering
  \begin{subfigure}[t]{0.48\textwidth}
    \includegraphics[width=\textwidth]{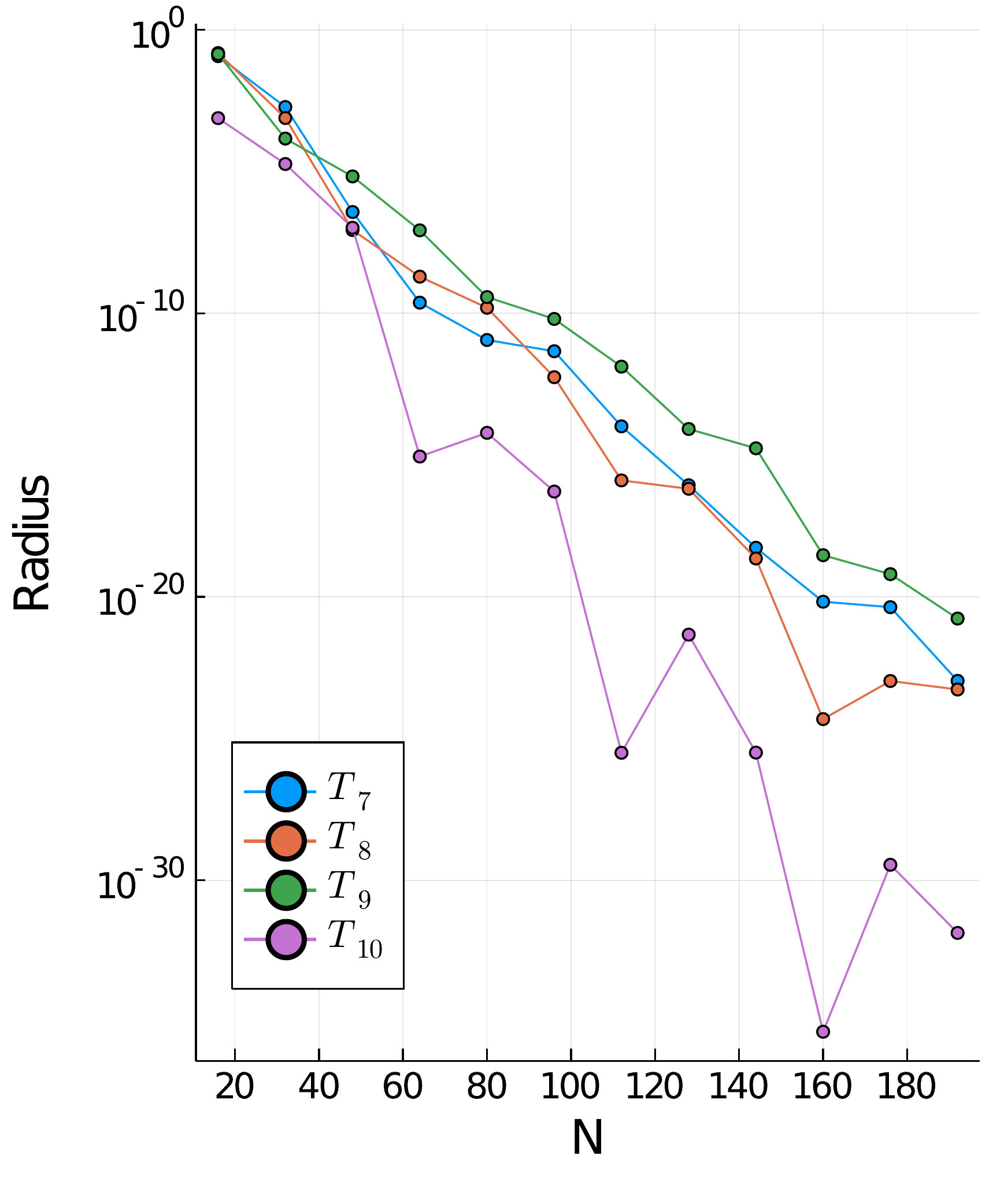}
    \caption{High precision computations of the enclosure}
    \label{fig:singular-triangles-high-precision-convergence}
  \end{subfigure}
  \hspace{0.02\textwidth}
  \begin{subfigure}[t]{0.48\textwidth}
    \includegraphics[width=\textwidth]{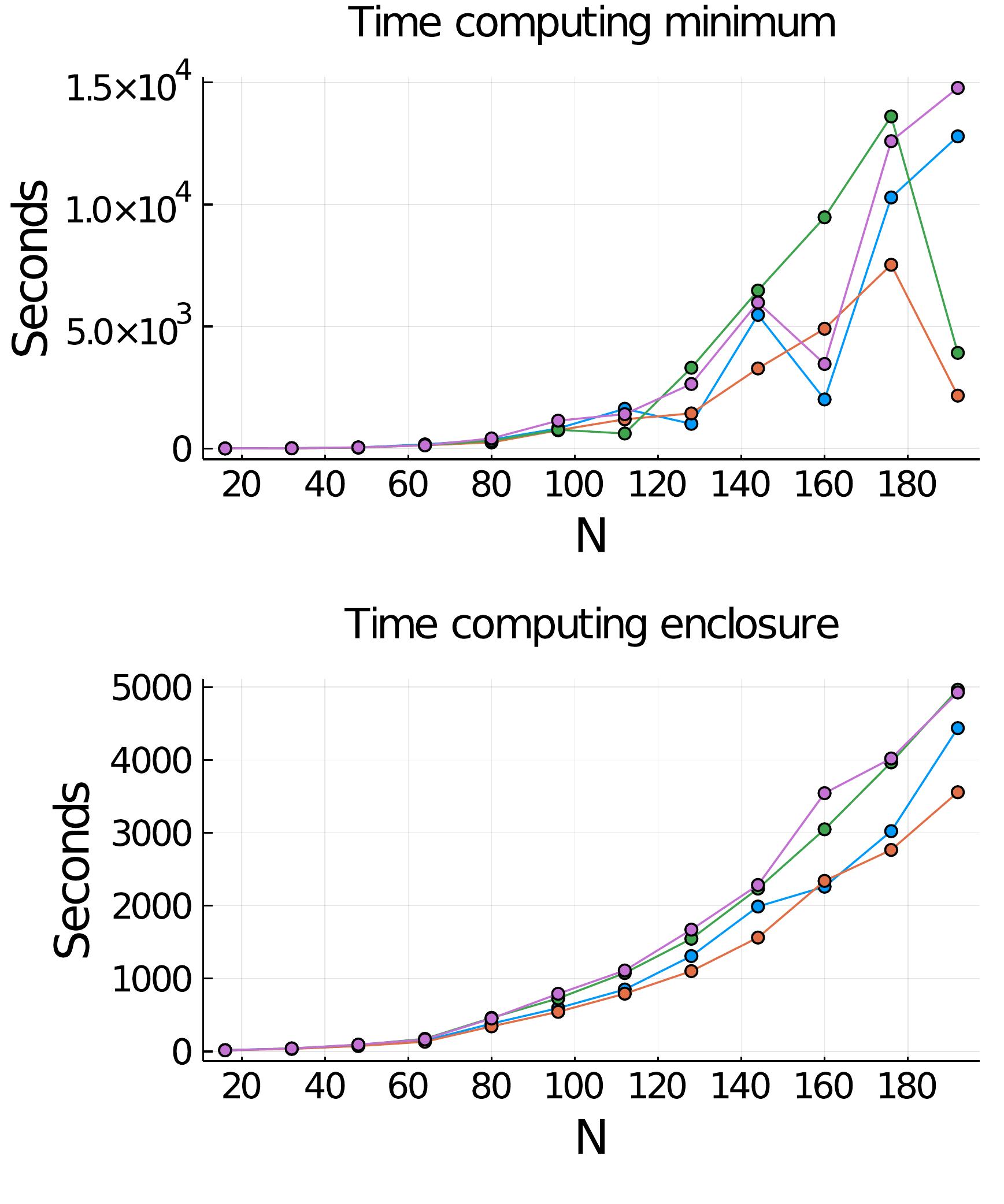}
    \caption{Time for computing the minimum of $\sigma(\lambda)$ as
      well as the enclosure.}
    \label{fig:singular-triangles-high-precision-time}
  \end{subfigure}
  \caption{High precision computations for the triangles in
    Table~\ref{table:singular-triangles}.}
  \label{fig:singular-triangles-high-precision}
\end{figure}

High precision results for the same triangles are presented in
Figure~\ref{fig:singular-triangles-high-precision-convergence} and
timings in Figure~\ref{fig:singular-triangles-high-precision-time}.
The situation is very similar to that of the regular triangles in
Figure~\ref{fig:regular-triangles-high-precision-convergence} and
\ref{fig:regular-triangles-high-precision-time} except that the rate
of convergence is lower. The computation of the norm takes up a larger
part of the time than for the regular triangles but is still dominated
by the computation of the minimum for high values of $N$.

\subsection{Combinatorial application: denominators of asymptotic
exponents}\label{sec:denominators}
Our motivation in this study is to obtain lower bounds on the
denominators the asymptotic exponents $\alpha$ in
Equation~\eqref{eq:asymptfn} would have if they were rational numbers.

This is achieved by the computation of a continued fraction expansion,
a routine technique in experimental mathematics. First, the
exponent~$\alpha$ is computed using interval arithmetic from the
enclosure of the fundamental eigenvalue. Next, the regular continued
fraction is computed using interval arithmetic and the computation is
stopped at the first time a partial quotient cannot be guaranteed. The
sequence of integers thus obtained is used to compute exactly the
first convergents $(P_n/Q_n)$ of the continued fraction. Unless the
number is thus detected to be rational, the last $Q_n$ is a lower
bound on the actual denominator. The computed lower bounds are listed
in Table~\ref{table:triangles-high-precision-denominator} in the
Appendix.

\subsection{3D Kreweras walks}\label{sec:3dkreweras}
In dimension~2, the Kreweras walks are walks with step set
$\{(-1,0), (0,-1),(1,1)\}$ whose study started with a 100-page
article by Kreweras~\cite{Kreweras1965} showing in particular that the
generating functions of interest are algebraic (and therefore solutions
of a linear differential equation), but the proof is far from trivial.
In dimension~2, all step sets with small steps (where each coordinate
has absolute value at most~1) having generating functions that are
solutions of a linear differential equation are also step sets for
which an associated group of the walk is
finite~\cite{BostanRaschelSalvy2014}. This is just an observation
obtained by looking at all possible cases. It is natural to wonder
whether a deeper relation between this group and the generating
functions could explain this property.

In dimension~3, the 3D Kreweras walks are defined as the natural
generalization of the 2D case, with step set
$\{(-1,0,0),(0,-1,0),(0,0,-1),(1,1,1)\}$, for which the analogous
group is again finite. That step set leads to the study of the
spherical triangle with angles $(2\pi/3,2\pi/3,2\pi/3)$ and its
fundamental eigenvalue~\cite{BogoselPerrollazRaschelTrotignon2020}.
The history of the knowledge on this eigenvalue, following
Bogosel
et al.~\cite{BogoselPerrollazRaschelTrotignon2020} and a private
communication of Bostan, is as follows:
\begin{itemize}
  \item[] $[5.15,5.16]$ in 2008 by Costabel (unpublished?);
  \item[] 5.158968860560663 in 2009 by Ratzkin and Treibergs~\cite{RatzkinTreibergs2009};
  \item[] 5.1606 in 2013 by Balakrishna~\cite{Balakrishna2013};
  \item[] 5.159145642466 in 2015 by Guttmann (unpublished);
  \item[] 5.1591452 in 2016 by Bacher, Kauers and Yatchak~\cite{BacherKauersYatchak2016};
  \item[] 5.159145642470 in 2020 by Bogosel, Perrollaz, Raschel and
  Trotignon~\cite{BogoselPerrollazRaschelTrotignon2020}.
\end{itemize}
We can now certify that the eigenvalue is actually
\begin{multline*}
  \lambda = 5.1591456424665417112216748625993501893151700566462081663\\
  0858031086922413365742186774243415327168103656498\dots
\end{multline*}
The corresponding exponent $\alpha = -1-\sqrt{\lambda+1/4}$ is
\begin{multline*}
  \alpha = -3.325757004174456250974540734758388852786843862030738206\\
  09206024964659686065647234082158565813950933996592\dots
\end{multline*}
with continued fraction
\begin{multline*}
  [-4; 1, 2, 14, 3, 100, 12, 102, 1, 5, 1, 2, 7, 6, 1, 11, 1, 6, 4, 1,\\
  8, 3, 3, 1, 1, 44, 8, 3, 1, 3, 5, 1, 1, 2, 1, 2, 1, 4, 1, 1, 1, 6,\\
  4, 1, 2, 1, 3, 2, 1, 15, 1, 17, 1, 2, 1, 2, 1, 1, 5, 1, 2, 2, 13, 1,\\
  3, 15, 2, 1, 2, 1, 6, 6, 2, 1, 1, 1, 1, 2, 3, 1, 1, 19, 5, 1, 4, 2,\\
  7, 1, 1, 5, 1, 23, 195, 1, 1, 3, 1, 1, 3, 1, 1, 1, 1, 1, 9, 2, \dots].
\end{multline*}
This implies that if it is a rational number, its denominator
must be at
least
\[9571644798056984399060418592860369800792627450626933>10^{51}.\]

The computation to high accuracy exploits the symmetry of the
triangle. We use expansions at all vertices together with one from
the interior, but compared to Section~\ref{sec:singular-triangles} we
make full use of the symmetry of the domain. All of the vertices are
symmetric and hence we expect the coefficients for their expansions to
be the same. In addition, only every second term will appear. For the
interior expansion we get a six-fold symmetry and only every sixth
term appears. This gives us the expansion
\begin{multline*}
  u(\theta, \phi) = \sum_{l = 1}^{N_{2}}b_{k}\cos(3(l-1)\alpha\phi_{\text{int}})
                    \mathsf{P}_{\nu}^{-3(l-1)\alpha}(\cos(\theta_{\text{int}}))\\
                  + \sum_{k = 1}^{N_{1}}c_{k}
                    \left(\sin((2(k - 1) + 1)\alpha\phi_{1})
                    \mathsf{P}_{\nu}^{-(2(k - 1) + 1)\alpha}(\cos
                    (\theta_{1}))\right.\\
                    \qquad\qquad+ \sin((2(k - 1) + 1)\alpha\phi_{2})
                    \mathsf{P}_
                    {\nu}^{-(2(k - 1) + 1)\alpha}(\cos(\theta_{2}))\\
                    +\left. \sin((2(k - 1) + 1)\alpha\phi_{3})
                    \mathsf{P}_
                    {\nu}^{-(2(k - 1) + 1)\alpha}(\cos(\theta_
                    {3}))\right).
\end{multline*}
One benefit with our method is that we do not have to prove that the
above expansion satisfies the required symmetries, we just see a
better convergence if it does. The only property of the expansion that
we use when computing the enclosure is that it behaves in the same
way on all
three boundaries. It is therefore sufficient to bound the maximum
on only one of them.

\section*{Acknowledgements}
We thank Nick Trefethen for several helpful discussions related to the
Method of Particular Solutions and Gerard Orriols for discussions
about how to handle the singular triangles. We are also thankful to
Kilian~Raschel and his co-authors who kindly shared preprints and
suggestions. BS was supported in part by De Rerum Natura
ANR-19-CE40-0018.

\begin{table}
  \centering
  \begin{tabular}{|c|p{0.9\textwidth}|}
    \hline
              & Eigenvalue \\ \hline
    $T_{1}$ & \seqsplit{12.40005165284337790528605341289663672073595731895}\\ \hline
    $T_{2}$ & \seqsplit{13.744355213213231835401121592138020782806650259631874894136332068957983025438961921160}\\ \hline
    $T_{3}$ & \seqsplit{20.571973537984730556625842153297}\\ \hline
    $T_{4}$ & \seqsplit{21.30940763019044525895348144123051777833684257714671661311314241820623854704023394191230205956761157788382983670637759893972691694122541330093667358027491678658694284070553504990811731549297257589768013675637}\\ \hline
    $T_{5}$ & \seqsplit{24.4569137962991116944804381447726828996079591315663692293441391578879515149}\\ \hline
    $T_{6}$ & \seqsplit{49.109945263284609919670343151508268353698425615333956068479546500637275248339988486176558994445206617439284515387218370698834970763269465605779603204345057}\\ \hline
    $T_{7}$ & \seqsplit{4.2617347552939870857522}\\ \hline
    $T_{8}$ & \seqsplit{5.15914564246654171122167486259935018931517005664620816630858031086922413365742186774243415327168103656498}\\ \hline
    $T_{9}$ & \seqsplit{6.24174833072633423680}\\ \hline
    $T_{10}$ & \seqsplit{6.77710805459830095738567415001383748}\\ \hline
  \end{tabular}
  \caption{Correctly rounded eigenvalues for the triangles from
    Table~\ref{table:regular-triangles} and
    \ref{table:singular-triangles}.}
  \label{table:triangles-high-precision-results}
\end{table}

\renewcommand{\arraystretch}{2}
\begin{table}
  \centering
  \begin{tabular}{|c|p{0.9\textwidth}|}
    \hline
              & Lower bound on denominator \\ \hline
    $T_{1}$ & \seqsplit{465867258515962084358692}$\ >10^{23}$\\ \hline
    $T_{2}$ & \seqsplit{48134549993161040519120418784541049178022}$\ >10^{40}$\\ \hline
    $T_{3}$ & \seqsplit{590595775643963}$\ >10^{14}$\\ \hline
    $T_{4}$ & \seqsplit{10073225747318443303256795082181812108714524296082956747116043882446421582505743104713297988178901824645}$\ >10^{103}$\\ \hline
    $T_{5}$ & \seqsplit{10653865792211960990143189115972047342}$\ >10^{37}$\\ \hline
    $T_{6}$ & \seqsplit{20195981375371512441648107892888975177795954858610023334567519476735622670369}$\ >10^{76}$\\ \hline
    $T_{7}$ & \seqsplit{76966517564}$\ >10^{10}$\\ \hline
    $T_{8}$ & \seqsplit{9571644798056984399060418592860369800792627450626933}$\ >10^{51}$\\ \hline
    $T_{9}$ & \seqsplit{4454060404}$\ >10^{9}$\\ \hline
    $T_{10}$ & \seqsplit{197533395012500053}$\ >10^{17}$\\ \hline
  \end{tabular}
  \caption{Lower bound on the denominator of the asymptotic exponent
    for the triangles from Table~\ref{table:regular-triangles} and
    \ref{table:singular-triangles}}
  \label{table:triangles-high-precision-denominator}
\end{table}

\end{document}